\documentclass[oneside,american]{amsart}
\usepackage[T1]{fontenc}
\usepackage{textcomp}
\usepackage[utf8]{inputenc}
\usepackage{mathrsfs}
\usepackage{mathtools}
\usepackage{amstext}
\usepackage{amsthm}
\usepackage{amssymb}
\usepackage{geometry}
\usepackage{esint}

\makeatletter
\numberwithin{equation}{section}
\numberwithin{figure}{section}
\theoremstyle{plain}
\newtheorem{thm}{\protect\theoremname}[section]
\theoremstyle{plain}
\newtheorem{lem}[thm]{\protect\lemmaname}
\ifx\proof\undefined
\newenvironment{proof}[1][\protect\proofname]{\par
	\normalfont\topsep6\p@\@plus6\p@\relax
	\trivlist
	\itemindent\parindent
	\item[\hskip\labelsep\scshape #1]\ignorespaces
}{%
	\endtrivlist\@endpefalse
}
\providecommand{\proofname}{Proof}
\fi
\theoremstyle{plain}
\newtheorem{prop}[thm]{\protect\propositionname}
\theoremstyle{plain}
\newtheorem{cor}[thm]{\protect\corollaryname}
\theoremstyle{remark}
\newtheorem{rem}[thm]{\protect\remarkname}

\makeatother

\usepackage{babel}
\providecommand{\corollaryname}{Corollary}
\providecommand{\lemmaname}{Lemma}
\providecommand{\propositionname}{Proposition}
\providecommand{\remarkname}{Remark}
\providecommand{\theoremname}{Theorem}

\begin{document}
\title[Quantum Variance and Benjamini-Schramm Convergence]{Quantum Variance of Holomorphic Modular Forms and Benjamini-Schramm
Convergence}
\author{Kemin Qiu}
\address{Magdalen College, Oxford, OX1 4AU, United Kingdom}
\email{kemin.qiu@magd.ox.ac.uk}
\begin{abstract}
We study quantum variance for holomorphic modular forms of fixed weight
on sequences of compact hyperbolic surfaces converging to the hyperbolic
plane in the sense of Benjamini-Schramm. We first show that quantum
ergodicity need not hold for an arbitrary choice of orthonormal basis,
even for sequences with a uniform lower bound on the injectivity radius
and a uniform spectral gap. The counterexample is constructed using
a tower of coverings and arises from the identity contribution in
the automorphic kernel. 

We then show that this obstruction is exceptional from the point of
view of the choice of basis. For a Haar random orthonormal basis of
holomorphic modular forms, the quantum variance converges to zero
with probability tending to one uniformly over all observables with
bounded supremum norm. The proof combines estimates for the Hilbert-Schmidt
norm of a projected multiplication operator with the Haar expectation
of the quantum variance, concentration of measure on the special unitary
group and a finite-dimensional approximation of the space of observables.
\end{abstract}

\maketitle

\section{Introduction}

\subsection{Related work }

Quantum ergodicity concerns the asymptotic distribution of eigenfunctions
of the Laplace-Beltrami operator on compact Riemannian manifolds whose
geodesic flow is ergodic. The classical theorem of Šnirel'man \cite{key-5},
Colin de Verdière \cite{key-6} and Zelditch \cite{key-7} asserts
that on a compact negatively curved manifold, as one takes Laplace
eigenfunctions with eigenvalues tending to infinity, almost all of
them become equidistributed. In particular, for compact hyperbolic
surfaces, Zelditch proved that as eigenvalue tends to infinity, the
mass of most eigenfunctions converges to the uniform measure.

A different asymptotic regime was studied by Le Masson and Sahlsten
\cite{key-4}, who proved a quantum ergodicity theorem for eigenfunctions
in a fixed spectral window along sequences of compact hyperbolic surfaces
converging to the hyperbolic plane in the sense of Benjamini and Schramm.
Under this geometric convergence, together with a uniform spectral
gap and injectivity radius lower bound, they showed that the quantum
variance of eigenfunctions in a fixed window tends to zero. Their
proof built on a wave propagation approach developed in earlier work
of Brooks, Le Masson and Lindenstrauss \cite{key-14} on discrete
graphs. A key step is to bound the quantum variance by the Hilbert-Schmidt
norm of an observable which is independent of the choice of eigenbasis.
This perspective inspires the first part of our work, where we similarly
upper bound the quantum variance for holomorphic modular forms in
terms of Hilbert-Schmidt norms of an operator related to the observables.

In the holomorphic setting, Rudnick and Sarnak \cite{key-8} formulated
an analogue of quantum unique ergodicity for Hecke eigenforms. For
the modular surface, the corresponding mass equidistribution theorem
in the weight aspect was proved by Holowinsky and Soundararajan \cite{key-9}.
Related equidistribution results in the level aspect were obtained
by Nelson \cite{key-10}. 

Random orthonormal bases provide another framework for quantum ergodicity.
Zelditch proved that a Haar-random orthonormal basis of spherical
harmonics on $S^{2}$ is quantum ergodic \cite{key-11}. He later
generalized this to a broad framework for sequences of high-dimensional
Hilbert spaces. He showed that a random orthonormal basis obtained
by applying a Haar-distributed unitary matrix to a fixed reference
basis, is almost surely quantum ergodic. A key ingredient in this
framework is the computation of moments of the quantum variance with
respect to the Haar measure on the unitary group, which uses the Weingarten
formula. Our Proposition \ref{prop:general meanQV} follows the same
approach in computing the expected quantum variance.

\subsection{Main results}

Let $\mathbb{D}$ be the Poincaré disk equipped with the hyperbolic
metric, and $X=\Gamma\backslash\mathbb{D}$ be a compact hyperbolic
surface. For a fixed integer $k>2$, let $\mathcal{M}_{k}\left(\Gamma\right)$
be the finite-dimensional space of holomorphic modular forms of weight
$k$ of $\Gamma$ equipped with Petersson inner product (defined in
Section \ref{sec:Setup-and-Notation}). Let
\[
d_{k}=\mathrm{dim}\mathcal{M}_{k}\left(\Gamma\right).
\]
For an orthonormal basis $\left\{ f_{j}\right\} _{j=1}^{d_{k}}$ of
$\mathcal{M}_{k}\left(\Gamma\right)$, and for an observable $a\in L^{\infty}\left(X\right)$,
we define the quantum variance of $a$ with respect to this basis
by 
\begin{equation}
\frac{1}{d_{k}}\sum_{j=1}^{d_{k}}\left|\left\langle f_{j},af_{j}\right\rangle _{k}-\fint ad\mathrm{Vol}_{X}\right|^{2}.\label{eq:quantum_variance}
\end{equation}
Here $\left\langle \cdot,\cdot\right\rangle _{k}$ is the Petersson
inner product and $\fint$ means that we divide the integral by the
volume of $X$. Since replacing $a$ by 
\[
a-\fint ad\mathrm{Vol}_{X}
\]
does not change the quantum variance, we may restrict attention to
observables satisfying
\[
\int ad\mathrm{Vol}_{X}=0.
\]
In this case, we can write the quantum variance (\ref{eq:quantum_variance})
as 
\begin{equation}
\frac{1}{d_{k}}\sum_{j=1}^{d_{k}}\left|\left\langle f_{j},af_{j}\right\rangle _{k}\right|^{2}\label{eq:qv_aa_n=00003D0}
\end{equation}

We are interested in the behavior of quantum variance in the large
volume limit for uniformly bounded observables. If a sequence of observables
$\left(a_{n}\right)_{n\in\mathbb{N}}$, where $a_{n}\in L^{\infty}\left(X_{n}\right)$,
satisfies 
\[
\left\Vert a_{n}\right\Vert _{\infty}\leq C
\]
for some constant $C$, we may replace $a_{n}$ by $a_{n}/C$, which
scales the quantum variance by $C^{-2}$ and does not affect whether
it tends to zero. Thus, throughout the paper, we assume without loss
of generality that 
\[
\sup_{n}\left\Vert a_{n}\right\Vert _{\infty}\leq1.
\]
Following Le Masson and Sahlsten \cite{key-4}, we consider sequences
of compact hyperbolic surfaces $X_{n}=\Gamma_{n}\backslash\mathbb{D}$
that converge to $\mathbb{D}$ in the sense of Benjamini and Schramm.
We denote by $\mathrm{InjRad}_{X}(z)$ the injectivity radius at a
point $z\in X$ and write $\mathrm{InjRad}(X)=\inf_{z\in X}\mathrm{InjRad}_{X}(z)$.
We say a sequence of compact hyperbolic surfaces $X_{n}=\Gamma_{n}\backslash\mathbb{D}$
converges to $\mathbb{D}$ in the Benjamini-Schramm sense if for any
$R>0$, 
\[
\lim_{n\to\infty}\frac{\mathrm{Vol}\left(\left\{ z\in X_{n}:\mathrm{InjRad}_{X_{n}}(z)<R\right\} \right)}{\mathrm{Vol}\left(X_{n}\right)}=0.
\]
Le Masson and Sahlsten proved that, for Maass forms in a fixed spectral
window, the quantum variance tends to zero as the surfaces converge
to $\mathbb{D}$ in the Benjamini-Schramm sense. One may naturally
ask whether an analogous statement holds for holomorphic modular forms,
with weight kept fixed while the underlying surface varies. In this
paper, we find that this does not always hold.
\begin{thm}
\label{thm:bad basis}There exists a sequence of compact hyperbolic
surfaces $X_{n}=\Gamma_{n}\backslash\mathbb{D}$ converging to $\mathbb{D}$
in the sense of Benjamini and Schramm, with 
\[
\mathrm{InjRad}\left(X_{n}\right)\geq r_{0}>0,
\]
for all $n$, where $r_{0}$ is a constant. For each $n$, one can
choose an orthonormal basis
\[
\mathcal{B}_{n,k}=\left\{ f_{n,1},\ldots,f_{n,d_{n,k}}\right\} 
\]
of $\mathcal{M}_{k}\left(\Gamma_{n}\right)$ and an observable $a_{n}$
with 
\[
\left\Vert a_{n}\right\Vert _{\infty}\leq1,
\]
such that 
\[
\liminf_{n\to\infty}\frac{1}{d_{n,k}}\sum_{j=1}^{d_{n,k}}\left|\left\langle f_{n,j},a_{n}f_{n,j}\right\rangle -\fint a_{n}d\mathrm{Vol}_{X_{n}}\right|^{2}>0,
\]
where 
\[
d_{n,k}=\mathrm{dim}\mathcal{M}_{k}\left(\Gamma_{n}\right).
\]
\end{thm}
We construct this counterexample by using a tower of coverings in
Section \ref{sec:Construction-of-Counterexample}. Moreover, the tower
of coverings used to construct the counterexample can be chosen to
have a uniform spectral gap, as explained in Remark \ref{rem:uniform spectral gap}.
Thus, the counterexample persists even under the uniform gap assumption
of Le Masson and Sahlsten\cite{key-4}. The obstruction to the convergence
of the quantum variance is the identity contribution to the Hilbert-Schmidt
norm, as we explain in Section \ref{sec:The-Hilbert-Schmidt-Norm}.

Thus, the failure of quantum ergodicity exhibited in Theorem \ref{thm:bad basis}
depends on the choice of basis. Despite the existence of such bad
bases, we prove in the second part of the paper that most bases are
quantum ergodic. To state this precisely, fix a reference basis for
$\mathcal{M}_{k}\left(\Gamma_{n}\right)$,
\[
\mathcal{E}_{n,k}=\left\{ e_{n,1},\ldots e_{n,d_{n,k}}\right\} .
\]
For an observable $a_{n}\in L^{\infty}\left(X_{n}\right)$ and $U\in\mathrm{SU}\left(d_{n,k}\right)$,
write 
\[
\mathcal{V}_{n,k}\left(U,a_{n}\right)\coloneqq\frac{1}{d_{n,k}}\sum_{j=1}^{d_{n,k}}\left|\left\langle Ue_{n,j},a_{n}Ue_{n,j}\right\rangle -\fint a_{n}d\mathrm{Vol}_{X_{n}}\right|^{2}.
\]
Let $\mu_{n,k}$ denote the normalized Haar probability measure on
$\mathrm{SU}\left(d_{n,k}\right)$, where $\mathrm{SU}(d)$ is the
group of $d\times d$ unitary matrices with determinant $1$.
\begin{thm}
\label{thm:most coverge}Let $X_{n}=\Gamma_{n}\backslash\mathbb{D}$
be a sequence of compact hyperbolic surfaces converging to $\mathbb{D}$
in the sense of Benjamini and Schramm, and suppose that 
\[
\mathrm{InjRad}\left(X_{n}\right)\geq r_{0}>0,
\]
 for each $n$ in $\mathbb{N}$. Fix $k>2$, and let $\mu_{n,k}$
denote the normalized Haar probability measure on $\mathrm{SU}\left(d_{n,k}\right)$.
Then for every $\varepsilon>0$,
\[
\mu_{n,k}\left(\left\{ U\in\mathrm{SU}\left(d_{n,k}\right):\sup_{\left\Vert a_{n}\right\Vert _{\infty}\leq1}\mathcal{V}_{n,k}\left(U,a_{n}\right)\leq\varepsilon\right\} \right)\to1
\]
as $n\to\infty$.
\end{thm}
Theorem \ref{thm:most coverge} shows that, with Haar probability
tending to one, the quantum variance tends to $0$ over all observables
of sup norm at most one. In particular, the result holds simultaneously
for all such observables, rather than for a single fixed observable.

The proof of Theorem \ref{thm:most coverge} combines the expectation
of the quantum variance over Haar-random bases with the Gromov-Milman
concentration inequality and a finite net argument for the low-frequency
part of the spherical decomposition. The high-frequency contribution
is shown to be negligible, while the remaining geometric error tends
to zero under Benjamini-Schramm convergence.

\subsection{Organization of the paper }

The two results are proved using a common Hilbert-Schmidt framework,
but the identity contribution plays opposite roles in the deterministic
and random-basis settings.

In Section \ref{sec:Setup-and-Notation}, we introduce the geometric
and spectral background. We then reformulate the problem on $\Gamma\backslash G$
by introducing a unitary lift from $\Gamma\backslash\mathbb{D}$ to
$\Gamma\backslash G$ in Section \ref{sec:Unitary-Lift-to}. We then
introduce the holomorphic projection $Q_{k}$ in Section \ref{sec:The-Holomorphic-Projection},
and upper bound the quantum variance in terms of the Hilbert-Schmidt
norm of $Q_{k}AQ_{k}$. The resulting Hilbert-Schmidt norm is then
decomposed into the identity and non-identity contributions. The explicit
decomposition is done in Section \ref{sec:The-Hilbert-Schmidt-Norm}.
We then use the identity contribution to construct a sequence of bases
for which quantum ergodicity fails, proving Theorem \ref{thm:bad basis}.

We then turn to Haar-random bases. Section \ref{sec:Expectation-Over-Haar-Random}
computes the expected quantum variance, and Section \ref{sec:Concentration-Over-Haar-Random}
combines this estimate with concentration of measure and a finite-dimensional
low-frequency approximation to obtain uniform quantum ergodicity for
most bases, proving Theorem \ref{thm:most coverge}.

\subsection*{Acknowledgements}

I would like to thank Michael Magee for introducing me to this project
and for many helpful discussions and suggestions. I would also like
to thank Will Hide for his careful reading of the manuscript and helpful
comments.

\section{Setup and Notation}\label{sec:Setup-and-Notation}

\subsection{Hyperbolic surfaces}

Let 
\[
\mathbb{D}=\left\{ z\in\mathbb{C}:\left|z\right|<1\right\} 
\]
be the Poincaré disk equipped with the hyperbolic metric
\[
ds^{2}=\frac{4\left|dz\right|^{2}}{\left(1-\left|z\right|^{2}\right)^{2}}.
\]
The associated volume form is
\[
d\mathrm{Vol}_{\mathbb{D}}(z)=\frac{4dxdy}{\left(1-\left|z\right|^{2}\right)^{2}},
\]
where $z=x+iy$. We denote the hyperbolic distance between $z,w\in\mathbb{D}$
by $d_{\mathbb{D}}\left(z,w\right)$. Let 
\[
G=\mathrm{SU}\left(1,1\right)=\left\{ g=\left(\begin{array}{cc}
\alpha(g) & \beta(g)\\
\overline{\beta(g)} & \overline{\alpha(g)}
\end{array}\right):\left|\alpha(g)\right|^{2}-\left|\beta(g)\right|^{2}=1\right\} 
\]
act on $\mathbb{D}$ by Möbius transformations 
\[
g\cdot z=\frac{\alpha(g)z+\beta(g)}{\overline{\beta(g)}z+\overline{\alpha(g)}}.
\]
The action is transitive, and the stabilizer of the origin $0\in\mathbb{D}$
is 
\[
K=\left\{ \left(\begin{array}{cc}
e^{i\theta} & 0\\
0 & e^{-i\theta}
\end{array}\right):\theta\in\mathbb{R}\right\} .
\]
Thus we may identify 
\[
\mathbb{D}\cong G/K.
\]
For a torsion-free cocompact Fuchsian subgroup $\Gamma<G$, set 
\[
X=\Gamma\backslash\mathbb{D},\qquad Y=\Gamma\backslash G.
\]
The volume form $d\mathrm{Vol}_{\mathbb{D}}$ is $\Gamma$-invariant
and descends to a volume form on $X$, which we denote by $d\mathrm{Vol}_{X}$. 

For $z\in X$, define the injectivity radius at $z$ by 
\[
\mathrm{InjRad}_{X}(z)=\frac{1}{2}\min_{\gamma\in\Gamma-\left\{ \mathrm{id}\right\} }d_{\mathbb{D}}\left(\tilde{z},\gamma\tilde{z}\right),
\]
where $\tilde{z}\in\mathbb{D}$ is any lift of $z$. 

For some of the representation-theoretic computations below, it is
convenient to work in the upper half-plane model
\[
\mathbb{H}=\left\{ z=x+iy\in\mathbb{C}:y>0\right\} .
\]
The Cayley transform 
\[
\mathcal{C}:\mathbb{D}\to\mathbb{H},\qquad\mathcal{C}(z)=i\frac{1+z}{1-z}
\]
is an isometry when $\mathbb{H}$ is equipped with the hyperbolic
metric 
\[
ds^{2}=\frac{dx^{2}+dy^{2}}{y^{2}}.
\]
The corresponding volume form is 
\[
d\mathrm{Vol}_{\mathbb{H}}(z)=\frac{dxdy}{y^{2}}.
\]
Let 
\[
C=\left(\begin{array}{cc}
i & i\\
-1 & 1
\end{array}\right)\in\mathrm{GL}_{2}\left(\mathbb{C}\right)
\]
be the matrix representing the Cayley transform. Then conjugation
by $C$ defines an isomorphism of Lie groups
\[
\mathrm{SU}\left(1,1\right)\to\mathrm{SL}_{2}\left(\mathbb{R}\right),\qquad\text{by }g\mapsto CgC^{-1}.
\]
Under this identification, the actions on $\mathbb{D}$ and $\mathbb{H}$
satisfy
\[
\mathcal{C}\left(g\cdot z\right)=\left(CgC^{-1}\right)\cdot\mathcal{C}(z).
\]

\subsection{Holomorphic modular forms}

A holomorphic function $f$ on $\mathbb{D}$ is a holomorphic modular
form of weight $k$ of $\Gamma$ if 
\[
f(g\cdot z)=j(g,z)^{k}f(z),
\]
 for all $g\in\Gamma$, where 
\[
j(g,z)=\overline{\beta}z+\overline{\alpha},\qquad g=\left(\begin{array}{cc}
\alpha & \beta\\
\overline{\beta} & \overline{\alpha}
\end{array}\right).
\]
 For $k>2$, let $\mathcal{M}_{k}\left(\Gamma\right)$ denote the
space of holomorphic modular forms of weight $k$ for $\Gamma$. The
space is finite dimensional and is equipped with the Petersson inner
product 
\[
\left\langle f,g\right\rangle _{k}=\int_{\Gamma\backslash\mathbb{D}}\overline{f(z)}g(z)\left(1-\left|z\right|^{2}\right)^{k}d\mathrm{Vol}_{X}(z).
\]
The inner product is well-defined because the factor $\left(1-\left|z\right|^{2}\right)^{k}$
compensates the transformation of holomorphic modular forms. Its dimension
is given by the Riemann-Roch theorem and Gauss-Bonnet formula, for
$k>2$ ,
\begin{equation}
d_{k}=\frac{k-1}{4\pi}\mathrm{Vol}\left(X\right).\label{eq:d_n,k and Vol}
\end{equation}

\subsection{Spectral decomposition of $L^{2}(X)$ }

Since $X$ is compact, the Laplace-Beltrami operator $\Delta$ on
$X$ has purely discrete spectrum. We denote its eigenvalues, counted
with multiplicity, by 
\[
0=\lambda_{0}<\lambda_{1}\leq\lambda_{2}\leq\ldots\to\infty.
\]
For each distinct eigenvalue $\lambda$, let $E_{\lambda}$ denote
the corresponding eigenspace. Then
\[
L^{2}\left(X\right)=\widehat{\bigoplus}_{\lambda\in\mathrm{Spec}\left(\Delta\right)}E_{\lambda},
\]
where the sum is the orthogonal Hilbert space direct sum over the
distinct eigenvalues.

Under the identification of $L^{2}(X)$ with the right $K$-invariant
subspace of $L^{2}\left(Y\right)$, the eigenspaces with $\lambda_{j}\geq1/4$
correspond to spherical principal-series representations, those with
$0<\lambda_{j}<1/4$ correspond to complementary-series representations,
and the eigenspace $E_{\lambda_{0}}$ corresponds to the trivial representation.

\section{Unitary Lift to $\Gamma\backslash G$ }\label{sec:Unitary-Lift-to}

Let $L^{2}(\Gamma\backslash\mathbb{D},k)$ denote the Hilbert space
of measurable functions $f:\mathbb{D\to C}$ satisfying 
\[
f(\gamma\cdot z)=j(\gamma,z)^{k}f(z),\qquad\text{for }\gamma\in\Gamma,
\]
equipped with the Petersson inner product. Let $L^{2}\left(\Gamma\backslash G,k\right)$
be the space of measurable functions $\phi:G\to\mathbb{C}$ satisfying
\[
\phi\left(\gamma g\kappa_{\theta}\right)=e^{ik\theta}\phi(g),\qquad\text{for }\gamma\in\Gamma,g\in G\text{ and }\kappa_{\theta}\in K,
\]
and 
\[
\int_{\Gamma\backslash G}\left|\phi\left(g\right)\right|^{2}dg<\infty.
\]
We normalize Haar measure on $G$ so that, writing 
\begin{equation}
g=g_{z}\cdot\kappa_{\theta},\qquad g_{z}\cdot0=z,\label{eq:quotient coordinate}
\end{equation}
one has
\begin{equation}
dg=\frac{4dxdy}{\left(1-\left|z\right|^{2}\right)^{2}}\cdot\frac{d\theta}{2\pi},\label{eq:Haar-measure}
\end{equation}
Then, $L^{2}(\Gamma\backslash G,k)$ and $L^{2}(\Gamma\backslash\mathbb{D},k)$
are unitarily isomorphic. An explicit unitary isomorphism is given
by
\[
\sigma_{k}:L^{2}(\Gamma\backslash\mathbb{D},k)\to L^{2}(\Gamma\backslash G,k),
\]
 where
\[
\left(\sigma_{k}f\right)(g)=\overline{\alpha(g)}^{-k}f\left(g\cdot0\right),
\]
Its inverse is given by 
\[
\left(\sigma_{k}^{-1}\phi\right)(z)=j(g,0)^{k}\phi(g),
\]
where $g$ is any element satisfying $g\cdot0=z$. The map $\sigma_{k}$
gives a unitary identification between the weighted space $L^{2}\left(\Gamma\backslash\mathbb{D},k\right)$
and the $K$-type $k$ subspace $L^{2}\left(\Gamma\backslash G,k\right)$.
Upon restriction to the holomorphic subspace, this agrees with the
standard disk realization of the holomorphic discrete series; see
Knapp \cite[Chapter VI, §1]{key-1}. 

Then for each $\Gamma_{n},$ let 
\[
\psi_{n,j}=\sigma_{k}f_{n,j},
\]
and lift the observable to a right $K$-invariant function on $\Gamma_{n}\backslash G$
by

\begin{equation}
A_{n}(g)=a_{n}(g\cdot0).\label{eq:A_n pull back}
\end{equation}
Then 
\[
\sigma_{k}\left(a_{n}f_{n,j}\right)=A_{n}\psi_{n,j}.
\]
Since $\sigma_{k}$ is unitary 
\[
\left\langle f_{n,j},a_{n}f_{n,j}\right\rangle _{k}=\left\langle \psi_{n,j},A_{n}\psi_{n,j}\right\rangle _{L^{2}\left(\Gamma_{n}\backslash G\right)}.
\]
Hence, for mean-zero observables, the quantum variance can be written
as 
\begin{equation}
\frac{1}{d_{n,k}}\sum_{j=1}^{d_{n,k}}\left|\left\langle \psi_{n,j},A_{n}\psi_{n,j}\right\rangle \right|^{2}.\label{eq:Gamma=00005CGqv}
\end{equation}

\section{The Holomorphic Projection}\label{sec:The-Holomorphic-Projection}

Define
\[
q_{k}:G\to\mathbb{C},\qquad q_{k}(g)=\overline{\alpha(g)}^{-k},
\]
and, for $\phi\in L^{2}(G)$,
\[
\left(Q_{k}\phi\right)(x)=\frac{k-1}{4\pi}\int_{G}q_{k}\left(y^{-1}x\right)\phi(y)dy.
\]
The operator $Q_{k}$ is the group realization of the weighted Bergman
projection. In particular 
\[
Q_{k}^{2}=Q_{k},\qquad Q_{k}^{*}=Q_{k}.
\]
Since $Q_{k}$ is defined by right convolution with a kernel on $G$,
the same operator acts on $L^{2}\left(\Gamma_{n}\backslash G\right)$
for every $n$. When regarded as an operator on $L^{2}\left(\Gamma_{n}\backslash G\right)$,
its image is precisely $\sigma_{k}\left(\mathcal{M}_{k}\left(\Gamma_{n}\right)\right)$.
Moreover, $Q_{k}$ acts as the identity on this subspace. As $Q_{k}$
is self-adjoint, we can write (\ref{eq:Gamma=00005CGqv}) as
\begin{equation}
\frac{1}{d_{n,k}}\sum_{j=1}^{d_{n,k}}\left|\left\langle \psi_{n,j},Q_{k}A_{n}Q_{k}\psi_{n,j}\right\rangle \right|^{2}.\label{eq:Q_kAQ_k qv}
\end{equation}
We can bound (\ref{eq:Q_kAQ_k qv}) by the Hilbert-Schmidt norm of
$Q_{k}A_{n}Q_{k}$ on $\Gamma_{n}\backslash G$, that is 
\begin{equation}
\frac{1}{d_{n,k}}\sum_{j=1}^{d_{n,k}}\left|\left\langle \psi_{n,j},Q_{k}A_{n}Q_{k}\psi_{n,j}\right\rangle \right|^{2}\leq\frac{1}{d_{n,k}}\left\Vert Q_{k}A_{n}Q_{k}\right\Vert _{\mathrm{HS}}^{2}.\label{eq:HS bound for qv}
\end{equation}

\section{The Hilbert-Schmidt Norm}\label{sec:The-Hilbert-Schmidt-Norm}

The kernel of $Q_{k}$ on $\Gamma\backslash G$ is the automorphic
Bergman kernel 
\[
B_{k}^{\Gamma}(x,y)=\frac{k-1}{4\pi}\sum_{\gamma\in\Gamma}q_{k}\left(y^{-1}\gamma x\right).
\]
Suppose that $A$ is the pullback of a function on $\Gamma\backslash\mathbb{D}$
under the quotient map $\Gamma\backslash G\to\Gamma\backslash\mathbb{D}$.
Then the kernel of $Q_{k}AQ_{k}$ on $\Gamma\backslash G$ can be
expressed as 
\[
K_{A}(x,y)=\int_{\Gamma\backslash G}B_{k}^{\Gamma}(x,z)A(z)B_{k}^{\Gamma}(z,y)dz.
\]
Then its Hilbert-Schmidt norm on $\Gamma\backslash G$ can be expressed
as 
\begin{align*}
\left\Vert Q_{k}AQ_{k}\right\Vert _{\mathrm{HS}}^{2} & =\int_{\Gamma\backslash G}\int_{\Gamma\backslash G}\left|\int_{\Gamma\backslash G}B_{k}^{\Gamma}(x,z)A(z)B_{k}^{\Gamma}(z,y)dz\right|^{2}dxdy\\
 & =\int_{\left(\Gamma\backslash G\right)^{4}}B_{k}^{\Gamma}(x,z)A(z)B_{k}^{\Gamma}(z,y)B_{k}^{\Gamma}(y,z')\overline{A(z')}B_{k}^{\Gamma}(z',x)dz'dzdxdy\\
 & =\int_{\Gamma\backslash G\times\Gamma\backslash G}A(z)\overline{A(z')}\left|B_{k}^{\Gamma}(z,z')\right|^{2}dz'dz,
\end{align*}
where the last equality is due to $Q_{k}^{2}=Q_{k}.$

We can expand $\left|B_{k}^{\Gamma}(x,y)\right|^{2}$ as 
\[
\left|B_{k}^{\Gamma}(x,y)\right|^{2}=\left(\frac{k-1}{4\pi}\right)^{2}\sum_{\delta,\eta\in\Gamma}q_{k}\left(y^{-1}\delta x\right)\overline{q_{k}\left(y^{-1}\eta x\right)},
\]
Set $\gamma=\delta\eta^{-1}$, and make the change of variable $x\mapsto\eta x$.
Since $A$ is left $\Gamma$-invariant, unfolding the integration
with respect $x$ gives us 
\[
\left\Vert Q_{k}AQ_{k}\right\Vert _{\mathrm{HS}}^{2}=\left(\frac{k-1}{4\pi}\right)^{2}\int_{\mathcal{F}}\int_{G}\sum_{\gamma\in\Gamma}A(x)\overline{A(y)}q_{k}\left(y^{-1}\gamma x\right)\overline{q_{k}\left(y^{-1}x\right)}dxdy,
\]
where $\mathcal{F}$ is the fundamental domain of $\Gamma$ on $G$.
Then 
\begin{equation}
\left\Vert Q_{k}AQ_{k}\right\Vert _{\mathrm{HS}}^{2}=\left(\frac{k-1}{4\pi}\right)^{2}\left(I_{e}+\sum_{\gamma\not=e}I_{\gamma}\right),\label{eq:HS=00003Di+n}
\end{equation}
where 
\begin{equation}
I_{\gamma}=\int_{\mathcal{F}}\int_{G}A(x)\overline{A(y)}q_{k}\left(y^{-1}\gamma x\right)\overline{q_{k}(y^{-1}x)}dxdy.\label{eq:I_=00005Cgamma}
\end{equation}

\subsection{The Non-identity Contribution}

We first estimate the contribution of non-identity elements of $\Gamma$.
It is clear from (\ref{eq:I_=00005Cgamma}) that
\begin{equation}
\left|\sum_{\gamma\not=e}I_{\gamma}\right|\leq\left\Vert A\right\Vert _{\infty}^{2}\sum_{\gamma\not=e}\int_{\mathcal{F}}\int_{G}\left|q_{k}\left(y^{-1}\gamma x\right)q_{k}(x^{-1}y)\right|dxdy\label{eq:bound 1}
\end{equation}

\begin{lem}
\label{lem:bound 2}Let 
\[
r_{\gamma}(y)=d_{\mathbb{D}}\left(y\cdot0,\gamma y\cdot0\right).
\]
Then,
\[
\int_{G}\left|q_{k}(y^{-1}\gamma x)q_{k}(x^{-1}y)\right|dx\lesssim_{k}\left(1+r_{\gamma}(y)\right)e^{-kr_{\gamma}(y)/2}.
\]
\end{lem}
\begin{proof}
Since $g\cdot0=\beta(g)/\overline{\alpha(g)}$ and $\left|\alpha\left(g\right)\right|^{2}-\left|\beta\left(g\right)\right|^{2}=1$,
we have 
\[
\left|\alpha(g)\right|=\left(1-\left|g\cdot0\right|^{2}\right)^{-1/2}.
\]
Hence 
\[
\left|q_{k}(g)\right|=\left(\cosh\frac{d_{\mathbb{D}}(0,g\cdot0)}{2}\right)^{-k}.
\]
As $\cosh t\geq\frac{e^{t}}{2}$, for all $t\in\mathbb{R}$, 
\begin{equation}
\left|q_{k}\left(g\right)\right|\leq2^{k}e^{-\frac{k}{2}d_{\mathbb{D}}\left(0,g\cdot0\right)}\label{eq:q_k upper bound}
\end{equation}
This implies 
\[
\left|q_{k}(y^{-1}\gamma x)q_{k}(x^{-1}y)\right|\leq4^{k}\exp\left[-\frac{k}{2}\left(d_{\mathbb{D}}\left(z,\gamma^{-1}y_{0}\right)+d_{\mathbb{D}}\left(z,y_{0}\right)\right)\right]
\]
where $y_{0}=y\cdot0$ and $z=x\cdot0$ is the $\mathbb{D}$-coordinate
of $x$ in (\ref{eq:quotient coordinate}). Then,
\begin{equation}
\int_{G}\left|q_{k}(y^{-1}\gamma x)q_{k}(x^{-1}y)\right|dx\leq\int_{\mathbb{D}}4^{k}\exp\left[-\frac{k}{2}\left(d_{\mathbb{D}}\left(z,\gamma^{-1}y_{0}\right)+d_{\mathbb{D}}(z,y_{0}\right)\right]d\mathrm{Vol}_{\mathbb{D}}(z).\label{eq:coshvofq-k}
\end{equation}
Let $\ell_{\gamma}$ be the complete geodesic passing through $y_{0}$
and $\gamma^{-1}y_{0}$, parametrized by arc length so that 
\[
\ell_{\gamma}(0)=y_{0},\qquad\ell_{\gamma}\left(r_{\gamma}(y)\right)=\gamma^{-1}y_{0}.
\]
We use Fermi coordinates $(s,t)$ with respect to $\ell_{\gamma}$,
where $s$is the signed arc length coordinate of orthogonal projection
of $z$ onto $\ell_{\gamma}$, measured from $y_{0}$, and $t$ is
the signed perpendicular distance from $z$ to $\ell_{\gamma}$. In
these coordinates, 
\[
d\mathrm{Vol_{\mathbb{D}}}(z)=\cosh tdsdt.
\]
By the hyperbolic right-angle formula 
\[
\cosh d_{\mathbb{D}}(z,y_{0})=\cosh s\cosh t,\ \cosh d_{\mathbb{D}}\left(z,\gamma^{-1}y_{0}\right)=\cosh\left(r_{\gamma}(y)-s\right)\cosh t.
\]
This implies that 
\[
d_{\mathbb{D}}\left(z,y_{0}\right)+d_{\mathbb{D}}\left(z,\gamma^{-1}y_{0}\right)\geq\left|s\right|+\left|r_{\gamma}(y)-s\right|+2\left|t\right|-C,
\]
for some constant $C$. Then for $s\in\left[0,r_{\gamma}(y)\right]$,
\[
d_{\mathbb{D}}\left(z,y_{0}\right)+d_{\mathbb{D}}\left(z,\gamma^{-1}y_{0}\right)\geq r_{\gamma}(y)+2\left|t\right|-C.
\]
For $s>r_{\gamma}(y)$ , 
\[
d_{\mathbb{D}}(z,y_{0})+d_{\mathbb{D}}\left(z,\gamma^{-1}y_{0}\right)\geq2s-r_{\gamma}(y)+2\left|t\right|-C=r_{\gamma}(y)+2\left(s-r_{\gamma}(y)\right)+2\left|t\right|-C.
\]
For $s<0$,
\[
d_{\mathbb{D}}\left(z,y_{0}\right)+d_{\mathbb{D}}\left(z,\gamma^{-1}y_{0}\right)\geq r_{\gamma}(y)-2s+2\left|t\right|-C.
\]
This implies, in Fermi coordinates, we can bound right hand side of
(\ref{eq:coshvofq-k}) by
\begin{multline*}
\int_{\mathbb{D}}\exp\left[-\frac{k}{2}\left(d_{\mathbb{D}}\left(z,\gamma^{-1}y_{0}\right)+d_{\mathbb{D}}(z,y_{0})\right)\right]dz\\
\lesssim\int_{-\infty}^{\infty}\exp\left(-\frac{k}{2}r_{\gamma}(y)-k\left|t\right|\right)\left(\int_{-\infty}^{0}e^{ks}ds+\int_{0}^{r_{\gamma}(y)}ds+\int_{r_{\gamma}(y)}^{\infty}e^{k\left(r_{\gamma}(y)-s\right)}ds\right)\cosh tdt.
\end{multline*}
This implies that 
\begin{align*}
\int_{G}\left|q_{k}(y^{-1}\gamma x)q_{k}(x^{-1}y)\right|dx & \lesssim\int_{-\infty}^{\infty}\exp\left(-\frac{k}{2}r_{\gamma}(y)-k\left|t\right|\right)\left(\frac{2}{k}+r_{\gamma}(y)\right)\cosh tdt\\
 & \lesssim_{k}\left(1+r_{\gamma}(y)\right)e^{-kr_{\gamma}(y)/2}.
\end{align*}
 
\end{proof}
\begin{prop}
\label{prop:bound-nonidentity}Let $X=\Gamma\backslash\mathbb{D}$
be a compact hyperbolic surface satisfying $\mathrm{InjRad}\left(X\right)\geq r_{0}>0$.
For any $R>0$ and $k>2$, we have 
\[
\left|\sum_{\gamma\not=e}I_{\gamma}\right|\leq C(k,r_{0})\left\Vert A\right\Vert _{\infty}^{2}\left[\sum_{m\geq2R}(m+2)e^{m-\frac{k}{2}m}\cdot\mathrm{Vol}(X)+\mathrm{Vol}\left(\left\{ z\in F:\mathrm{InjRad}_{X}(z)\leq R\right\} \right)\right],
\]
where $F$ is the fundamental domain of $\Gamma$ on $\mathbb{D}$,
$X=\Gamma\backslash\mathbb{D}$ and $C(k,r_{0})$ is a constant depending
only on $k$ and $r_{0}$.
\end{prop}
\begin{proof}
Let $p:G\to\mathbb{D}$ be the composition of the quotient map $G\to G/K$
and homeomorphism $G/K\to\mathbb{D}.$ Then we have 
\[
p^{-1}\left(F\right)=\mathcal{F},
\]
By using the coordinate in (\ref{eq:quotient coordinate}), we have
for a function $f$ that is right $K$-invariant 
\[
\int_{\mathcal{F}}f(g)dg=\int_{F}f\left(g_{z}\right)d\mathrm{Vol}_{\mathbb{D}}z.
\]
Then by (\ref{eq:bound 1}) and Lemma \ref{lem:bound 2}, 
\begin{align}
\left|\sum_{\gamma\not=e}I_{\gamma}\right| & \lesssim_{k}\sum_{\gamma\not=e}\left\Vert A\right\Vert _{\infty}^{2}\int_{\mathcal{F}}\left(1+r_{\gamma}(y)\right)e^{-kr_{\gamma}(y)/2}dy\label{eq:Bound 3}\\
 & =\sum_{\gamma\not=e}\left\Vert A\right\Vert _{\infty}^{2}\int_{F}\left(1+d_{\mathbb{D}}\left(z,\gamma z\right)\right)\exp\left(-\frac{k}{2}d_{\mathbb{D}}\left(z,\gamma z\right)\right)dz.\nonumber 
\end{align}
The last equality is because $r_{\gamma}(y)$ is right $K$-invariant.

Let 
\[
D(R)=\left\{ z\in F:\mathrm{InjRad}_{X}(z)>R\right\} .
\]
We split the geometric integral into contributions from $D\left(R\right)$
and $D\left(R\right)^{c}$.
\begin{multline*}
\sum_{\gamma\not=e}\int_{F}\left(1+d_{\mathbb{D}}\left(z,\gamma z\right)\right)\exp\left(-\frac{k}{2}d_{\mathbb{D}}\left(z,\gamma z\right)\right)dz=\int_{D(R)}\sum_{\gamma\not=e}\left(1+d_{\mathbb{D}}\left(z,\gamma z\right)\right)\exp\left(-\frac{k}{2}d_{\mathbb{D}}\left(z,\gamma z\right)\right)dz\\
+\int_{D(R)^{c}}\sum_{\gamma\not=e}\left(1+d_{\mathbb{D}}\left(z,\gamma z\right)\right)\exp\left(-\frac{k}{2}d_{\mathbb{D}}\left(z,\gamma z\right)\right)dz.
\end{multline*}

Let 
\[
N_{z}(m)=\#\left\{ \gamma\in\Gamma-\left\{ e\right\} :m\leq d_{\mathbb{D}}(z,\gamma z)<m+1\right\} .
\]
Let $r_{0}>0$ be a lower bound for the injectivity radius of $X$,
then for distinct $\gamma_{1},\gamma_{2}\in\Gamma,$
\[
d_{\mathbb{D}}\left(\gamma_{1}z,\gamma_{2}z\right)=d_{\mathbb{D}}\left(z,\gamma_{1}^{-1}\gamma_{2}z\right)\geq2r_{0}.
\]
Hence the balls $B\left(\gamma z,r_{0}/2\right)$ are pairwise disjoint.
For every $\gamma$ counted by $N_{z}\left(m\right)$, these balls
are contained in $B\left(z,m+1+\frac{r_{0}}{2}\right)$. Therefore
\[
N_{z}\left(m\right)\mathrm{Vol}\left(B\left(\frac{r_{0}}{2}\right)\right)\leq\mathrm{Vol}\left(B\left(m+1+\frac{r_{0}}{2}\right)\right).
\]
This implies that there exists a constant depending on $r_{0}$ such
that
\begin{equation}
N_{z}\left(m\right)\leq\frac{\cosh\left(m+1+\frac{r_{0}}{2}\right)-1}{\cosh\left(\frac{r_{0}}{2}\right)-1}\leq C_{1}\left(r_{0}\right)e^{m}.\label{eq:N_zm upper bound}
\end{equation}
Then for any $z$, we have 
\begin{align}
\sum_{\gamma\not=e}\left(1+d_{\mathbb{D}}\left(z,\gamma z\right)\right)\exp\left(-\frac{k}{2}d_{\mathbb{D}}\left(z,\gamma z\right)\right) & \leq\sum_{m=0}^{\infty}N_{z}(m)\left(m+2\right)e^{-km/2}.\label{eq:upperboundofnonidentityterm}
\end{align}
Since for $z\in D(R)$, $N_{z}(m)=0$ for $m<2R-1$, this implies
\begin{align*}
\sum_{\gamma\not=e}\left(1+d_{\mathbb{D}}\left(z,\gamma z\right)\right)\exp\left(-\frac{k}{2}d_{\mathbb{D}}\left(z,\gamma z\right)\right) & \leq\sum_{m\geq2R-1}C_{1}\left(r_{0}\right)\left(m+2\right)e^{m-km/2}\\
 & =C_{1}\left(r_{0}\right)e^{k/2-1}\sum_{m\geq2R}\left(m+1\right)e^{m-km/2}\\
 & \leq C_{1}\left(r_{0}\right)e^{k/2-1}\sum_{m\geq2R}\left(m+2\right)e^{m-km/2}
\end{align*}
where the series on right hand side converges for $k>2$. Then there
exists a constant $C_{2}\left(k,r_{0}\right)$ depending only on $k$
and $r_{0}$,
\begin{equation}
\int_{D(R)}\sum_{\gamma\not=e}\left(1+d_{\mathbb{D}}\left(z,\gamma z\right)\right)\exp\left(-\frac{k}{2}d_{\mathbb{D}}\left(z,\gamma z\right)\right)dz\leq C_{2}\left(k,r_{0}\right)\sum_{m\geq2R}(m+2)e^{m-km/2}\mathrm{Vol}(X).\label{eq:D(R) Bound}
\end{equation}
Similarly, there exists a constant $C_{3}(k,r_{0})$,
\begin{multline}
\int_{D(R)^{c}}\sum_{\gamma\not=e}\left(1+d_{\mathbb{D}}\left(z,\gamma z\right)\right)\exp\left(-\frac{k}{2}d_{\mathbb{D}}\left(z,\gamma z\right)\right)dz\\
\leq C_{3}\left(k,r_{0}\right)\left(\sum_{m=0}^{\infty}(m+2)e^{m-km/2}\right)\mathrm{Vol}\left(\left\{ z\in F:\mathrm{InjRad}_{X}(z)\leq R\right\} \right).\label{eq:D(R)^c bound}
\end{multline}
Combining (\ref{eq:D(R) Bound}) and (\ref{eq:D(R)^c bound}) completes
the proof.
\end{proof}
Then the following corollary follows directly from Proposition \ref{prop:bound-nonidentity}.
\begin{cor}
\label{cor:BS nonidentity term}Suppose that 
\[
X_{n}=\Gamma_{n}\backslash\mathbb{D}
\]
 converges to $\mathbb{D}$ in the sense of Benjamini-Schramm , and
assume that there exists $r_{0}>0$ such that 
\[
\mathrm{InjRad}\left(X_{n}\right)\geq r_{0}
\]
for all $n$. If $\left(A_{n}\right)_{n\in\mathbb{N}}$ is uniformly
bounded in $L^{\infty}$-norm, then 
\[
\frac{1}{\mathrm{Vol}\left(X_{n}\right)}\left|\sum_{\gamma\not=e}I_{n,\gamma}\right|\to0\qquad\text{as }n\to\infty.
\]
\end{cor}

\subsection{The Identity Contribution}

For the identity term,
\[
I_{e}=\int_{\mathcal{F}}\int_{G}A(x)\overline{A(y)}\left|q_{k}\left(y^{-1}x\right)\right|^{2}dxdy.
\]
Write 
\[
h_{k}(g)=\left|q_{k}(g)\right|^{2},
\]
and let $\rho(h_{k})$ be the operator acting on $L^{2}\left(\Gamma\backslash G\right)$
by 
\[
\left(\rho\left(h_{k}\right)f\right)(y)=\int_{G}h_{k}(x)\rho(x)f(y)dx=\int_{G}h_{k}(x)f(yx)dx.
\]
Then we can express the identity term as 
\begin{equation}
I_{e}=\left\langle A,\rho\left(h_{k}\right)A\right\rangle _{L^{2}\left(\Gamma\backslash G\right)}.\label{eq:I_e_1}
\end{equation}

Since $A$ comes from $L^{2}(X)$, we may use the spectral decomposition
from Section 2.3 and write $A$ as a sum of its spherical components.
Since $A$ is mean-zero, it is orthogonal to the trivial representation.
Hence
\begin{equation}
A=\sum_{\sigma\not\simeq\mathbf{1}}\sum_{a=1}^{m(\sigma,\Gamma)}A_{\sigma,a},\label{eq:decomp of A}
\end{equation}
where $A_{\sigma,a}$ denotes the component of $A$ in the $a$-th
copy of the spherical representation $\sigma$. 

Under the Cayley transform $G\simeq\mathrm{SL}_{2}\left(\mathbb{R}\right)$,
we use the realization of the spherical principal and complementary
series introduced by Repka \cite{key-2}. Let 
\[
M=\left\{ \pm I\right\} ,\qquad K=\mathrm{SO}(2),
\]
 and let $D$ and $N$ denote respectively the subgroups of positive
diagonal matrices and upper triangular unipotent matrices. Let 
\[
h_{t}=\left(\begin{array}{cc}
e^{t} & 0\\
0 & e^{-t}
\end{array}\right)\in D,\qquad n_{u}=\left(\begin{array}{cc}
1 & u\\
0 & 1
\end{array}\right)\in N.
\]
 Set 
\[
K^{+}=\left\{ \kappa_{\theta}=\left(\begin{array}{cc}
\cos\theta & -\sin\theta\\
\sin\theta & \cos\theta
\end{array}\right):0\leq\theta<\pi\right\} .
\]
Then every $g\in G$ can be uniquely written in the form 
\[
g=m(g)e^{H(g)}N(g)\kappa(g),
\]
 where $m(g)\in M$, $\kappa(g)\in K^{+}$, $N(g)\in N$ and $H(g)=\mathrm{diag}\left(t(g),-t(g)\right).$

Define 
\[
\mu\left(h_{t}\right)=e^{2t},\qquad\eta_{s}\left(mh_{t}n\right)=e^{st},
\]
 where we take the trivial character of $M$ here. Extend these functions
to $G$ by 
\[
\mu\left(g\right)=\mu\left(e^{H(g)}\right),\qquad\eta_{s}\left(g\right)=\eta_{s}\left(m(g)e^{H(g)}N(g)\right).
\]
Let 
\[
\pi_{s}=\mathrm{Ind}_{MDN}^{G}\eta_{s}.
\]
 The representation acts on $\mathcal{H}_{s}=L^{2}\left(K^{+}\right)$
by 
\begin{equation}
\left(\pi_{s}(g)f\right)(\kappa_{0})=\eta_{s}\left(\kappa_{0}g\right)\mu\left(\kappa_{0}g\right)^{1/2}f\left(\kappa(\kappa_{0}g)\right).\label{eq:pi s}
\end{equation}
For $s\in i\mathbb{R}$, $\pi_{s}$ is unitary and belongs to the
principal series. For $-1<s<0$, we equip the linear span of $\left\{ \phi_{2m}\right\} _{m\in\mathbb{Z}}$
with the inner product
\begin{equation}
\left\langle \phi_{2m},\phi_{2n}\right\rangle _{u,s}=\delta_{m,n}\begin{cases}
1, & m=0,\\
\prod_{j=1}^{\left|m\right|}\frac{2j-1-s}{2j-1+s} & m\not=0,
\end{cases}\label{eq:complementary series}
\end{equation}
where 
\[
\phi_{2m}\left(\kappa_{\theta}\right)=e^{2mi\theta},
\]
and $\left\{ \phi_{2m}\right\} _{m\in\mathbb{Z}}$ is an orthonormal
basis of $L^{2}\left(K^{+}\right)$ with respect to the usual $L^{2}$-inner
product. Let $\mathcal{H}_{s}^{u}$ be its Hilbert-space completion
with respect to $\left\langle \cdot,\cdot\right\rangle _{u,s}$. The
same action $\pi_{s}$ extends to a unitary representation on $\mathcal{H}_{s}^{u}$,
denoted $\pi_{s}^{u}$. This is the complementary series representation.
In both the principal and complementary series, the $K$-invariant
subspace is one dimensional and is spanned by 
\[
\phi_{0}\equiv1.
\]

\begin{lem}
\label{lem:rho(h_k)  fixes H^K}Let $\mathcal{H}_{s}$ be the spherical
principal or complementary series component of $L^{2}\left(\Gamma\backslash G\right)$.
Then 
\[
\rho\left(h_{k}\right)\mathcal{H}_{s}^{K}\subseteq\mathcal{H}_{s}^{K},
\]
 where $\mathcal{H}_{s}^{K}$ denotes the $K$-invariant subspace
of $\mathcal{H}_{s}$. In particular, since 
\[
\mathrm{dim}\mathcal{H}_{s}^{K}=1,
\]
 there exists $c_{k}(s)\in\mathbb{C}$ such that 
\[
\rho\left(h_{k}\right)v=c_{k}(s)v,\qquad v\in\mathcal{H}_{s}^{K}.
\]
\end{lem}
\begin{proof}
Let $v\in\mathcal{H}_{s}^{K}$. Since $\mathcal{H}_{s}$ is closed
under $G$-action, for each $g$ 
\[
h_{k}(g)\rho(g)v\in\mathcal{H}_{s}.
\]
As $\mathcal{H}_{s}$ is complete, 
\[
\rho\left(h_{k}\right)v=\int_{G}h_{k}(g)\rho(g)vdg\in\mathcal{H}_{s}.
\]
For $\kappa\in K$,
\begin{align*}
\rho\left(\kappa\right)\rho\left(h_{k}\right)v & =\int_{G}h_{k}(g)\rho\left(\kappa g\right)vdg\\
 & =\int_{G}h_{k}\left(\kappa^{-1}g\right)\rho(g)vdg\\
 & =\int_{G}h_{k}\left(g\right)\rho\left(g\right)vdg\\
 & =\rho\left(h_{k}\right)v.
\end{align*}
Hence it is $K$-invariant.
\end{proof}
Then by using the spectral decomposition of $A$ in (\ref{eq:decomp of A})
\[
\rho\left(h_{k}\right)A=\sum_{\sigma_{s}\not\simeq\mathbf{1}}\sum_{a=1}^{m(\sigma_{s},\Gamma)}\rho\left(h_{k}\right)A_{\sigma_{s},a}=\sum_{\sigma_{s}\not\simeq\mathbf{1}}\sum_{a=1}^{m(\sigma_{s},\Gamma)}c_{k}(s)A_{\sigma_{s},a},
\]
where $\sigma_{s}$ is isomorphic to the representation $\pi_{s}$
in the principal-series case and $\pi_{s}^{u}$ in the complementary-series
case. Then by (\ref{eq:I_e_1}) and orthogonality, 
\begin{equation}
I_{e}=\sum_{\sigma_{s}\not\simeq\mathbf{1}}\sum_{a=1}^{m\left(\sigma_{s},\Gamma\right)}c_{k}(s)\left\Vert A_{\sigma_{s},a}\right\Vert ^{2}.\label{eq:I_e_2}
\end{equation}

\begin{prop}
\label{prop:c_k calculation }Let $\sigma_{s}$ be the principal-series
representation with $s\in i\mathbb{R}$, or a complementary series
representation with $-1<s<0$. Then, 
\begin{equation}
c_{k}\left(s\right)=\int_{G}h_{k}(g)\left\langle e_{s},\sigma_{s}(g)e_{s}\right\rangle dg,\label{eq:c_k_G}
\end{equation}
where $e_{s}$ is the image of $1$ in $\mathcal{H}_{s}$ or $\mathcal{H}_{s}^{u}$
under the isomorphism. And this constant is given by 
\[
c_{k}(s)=4\pi\frac{\Gamma\left(k-\frac{1+s}{2}\right)\Gamma\left(k-\frac{1-s}{2}\right)}{\Gamma\left(k\right)^{2}}.
\]
\end{prop}
\begin{proof}
Since 
\[
c_{k}(s)=\left\langle e_{s},\rho\left(h_{k}\right)e_{s}\right\rangle ,
\]
the first part is clear. To calculate this we consider the Cartan
decomposition of $\mathrm{SL}_{2}\left(\mathbb{R}\right)$, which
is 
\[
\mathrm{SL_{2}}(\mathbb{R})=KDK.
\]
The coordinate of $g$ under this decomposition is 
\[
g=\kappa_{\phi}a_{r}\kappa_{\theta},\qquad a_{r}=\left(\begin{array}{cc}
e^{r/2} & 0\\
0 & e^{-r/2}
\end{array}\right),
\]
and the respective Haar measure under this coordinate is 
\[
dg=\sinh r\frac{d\theta}{2\pi}d\phi dr.
\]
Since $h_{k}$ is $K$-bi-invariant and $e_{s}$ is $K$-invariant,
we can express (\ref{eq:c_k_G}) as
\[
c_{k}(s)=2\pi\int_{0}^{\infty}\left(\cosh\frac{r}{2}\right)^{-2k}\left\langle e_{s},\sigma_{s}\left(a_{r}\right)e_{s}\right\rangle \sinh rdr.
\]
Using the unitary isomorphism, if $\sigma_{s}$ belongs to principal
series, then
\[
\left\langle e_{s},\sigma_{s}\left(a_{r}\right)e_{s}\right\rangle =\left\langle 1,\pi_{s}\left(a_{r}\right)1\right\rangle _{L^{2}}=\frac{1}{\pi}\int_{0}^{\pi}\left(\pi_{s}\left(a_{r}\right)1\right)\left(\kappa_{\theta}\right)d\theta.
\]
If $\sigma_{s}$ is a complementary series, we may expand $\pi_{s}\left(a_{r}\right)1$as
a sum of $\left\{ \phi_{2m}\right\} _{m\in\mathbb{Z}}$, 
\[
\pi_{s}\left(a_{r}\right)1=\sum_{m\in\mathbb{Z}}b_{2m}(r)\phi_{2m}.
\]
The orthogonality relation in (\ref{eq:complementary series}) implies
\[
\left\langle 1,\pi_{s}\left(a_{r}\right)1\right\rangle _{u,s}=b_{0}\left(r\right)=\left\langle 1,\pi_{s}\left(a_{r}\right)1\right\rangle _{L^{2}}.
\]
Therefore for all $s$, we may write 
\[
c_{k}(s)=2\pi\int_{0}^{\infty}\left(\cosh\frac{r}{2}\right)^{-2k}\left\langle 1,\pi_{s}\left(a_{r}\right)1\right\rangle _{L^{2}}\sinh rdr.
\]
By (\ref{eq:pi s}), we have 
\[
\left(\pi_{s}\left(a_{r}\right)1\right)\left(\kappa_{0}\right)=\mu\left(\kappa_{0}a_{r}\right)^{1/2}\eta_{s}\left(\kappa_{0}a_{r}\right).
\]
Using the Iwasawa decomposition of $\kappa_{\theta}a_{r}$, we obtain
\[
\left(\pi_{s}\left(a_{r}\right)1\right)\left(\kappa_{\theta}\right)=\left(e^{r}\sin^{2}\theta+e^{-r}\cos^{2}\theta\right)^{-\frac{1+s}{2}}=\left(\cosh r-\sinh r\cos2\theta\right)^{-\frac{1+s}{2}}.
\]
Then,
\begin{align*}
c_{k}(s) & =2\int_{0}^{\infty}\int_{0}^{\pi}\left(\cosh\frac{r}{2}\right)^{-2k}\left(\cosh r-\sinh r\cos2\theta\right)^{-\frac{1+s}{2}}\sinh rd\theta dr.
\end{align*}
Set 
\[
a=\frac{1+s}{2},\qquad u=\tanh^{2}\frac{r}{2}.
\]
Then 
\[
\left(\cosh\frac{r}{2}\right)^{-2k}=\left(1-u\right)^{k},
\]
 and 
\[
\cosh r-\sinh r\cos2\theta=\frac{1+u-2\sqrt{u}\cos2\theta}{1-u}.
\]
Moreover, 
\[
\sinh rdr=\frac{2du}{\left(1-u\right)^{2}}.
\]
Therefore
\[
c_{k}(s)=4\int_{0}^{1}\int_{0}^{\pi}\left(1-u\right)^{k+a-2}\left(1+u-2\sqrt{u}\cos2\theta\right)^{-a}d\theta du.
\]
Now
\[
1+u-2\sqrt{u}\cos2\theta=\left(1-\sqrt{u}e^{2i\theta}\right)\left(1-\sqrt{u}e^{-2i\theta}\right).
\]
Using the binomial expansion, 
\[
\left(1-z\right)^{-a}=\sum_{m=0}^{\infty}\frac{\left(a\right)_{m}}{m!}z^{m},\qquad\left|z\right|<1,
\]
 we obtain 
\[
\int_{0}^{\pi}\left(1+u-2\sqrt{u}\cos2\theta\right)^{-a}d\theta=\int_{0}^{\pi}\sum_{m,n\geq0}\frac{\left(a\right)_{m}\left(a\right)_{n}}{m!n!}u^{\left(m+n\right)/2}e^{2i\left(m-n\right)\theta}d\theta.
\]
Since 
\[
\int_{0}^{\pi}e^{2i\left(m-n\right)\theta}d\theta=\begin{cases}
\pi, & m=n,\\
0, & m\not=n,
\end{cases}
\]
it follows that 
\[
\int_{0}^{\pi}\left(1+u-2\sqrt{u}\cos2\theta\right)^{-a}d\theta=\pi\sum_{m=0}^{\infty}\frac{\left(a\right)_{m}^{2}}{\left(m!\right)^{2}}u^{m}.
\]
Hence 
\[
c_{k}(s)=4\pi\sum_{m=0}^{\infty}\frac{\left(a\right)_{m}^{2}}{\left(m!\right)^{2}}\int_{0}^{1}u^{m}\left(1-u\right)^{k+a-2}du.
\]
 By the beta-integral formula, 
\[
\int_{0}^{1}u^{m}\left(1-u\right)^{k+a-2}du=\frac{\Gamma\left(m+1\right)\Gamma\left(k+a-1\right)}{\Gamma\left(m+k+a\right)}.
\]
Thus
\begin{align*}
c_{k}(s) & =4\pi\frac{\Gamma\left(k+a-1\right)}{\Gamma\left(k+a\right)}\sum_{m=0}^{\infty}\frac{\left(a\right)_{m}^{2}}{m!\left(k+a\right)_{m}}\\
 & =4\pi\frac{\Gamma\left(k+a-1\right)\Gamma\left(k-a\right)}{\Gamma\left(k\right)^{2}},
\end{align*}
 where the last equality is due to Gauss's summation formula. Putting
\[
a=\frac{1+s}{2}
\]
gives us the result.
\end{proof}

\section{Construction of Counterexample}\label{sec:Construction-of-Counterexample}

Let $\Gamma_{0}<G$ be a torsion-free cocompact subgroup. Choose a
decreasing sequence of finite-index normal subgroups
\[
\Gamma_{0}\supset\Gamma_{1}\supset\Gamma_{2}\supset\ldots,\qquad\left[\Gamma_{0}:\Gamma_{n}\right]<\infty,\qquad\Gamma_{n}\vartriangleleft\Gamma_{0},
\]
such that 
\[
\bigcap_{n=1}^{\infty}\Gamma_{n}=\left\{ e\right\} .
\]
Since $\Gamma_{n}\subset\Gamma_{0}$, we have 
\[
\mathrm{InjRad}\left(X_{n}\right)\geq\mathrm{InjRad}\left(X_{0}\right)>0.
\]
Since $X_{0}$ is compact, there are only finitely many non-trivial
conjugacy classes in $\Gamma_{0}$ whose translation length is at
most $2R$. Let $\left\{ \gamma_{1},\ldots,\gamma_{N}\right\} $ be
representatives of these classes. Since 
\[
\bigcap_{n=1}^{\infty}\Gamma_{n}=\left\{ e\right\} ,
\]
for sufficiently large $n$, none of the elements $\gamma_{1},\ldots,\gamma_{N}$
belongs to $\Gamma_{n}$. As $\Gamma_{n}$ is normal in $\Gamma_{0}$,
no conjugate of any $\gamma_{i}$ belongs to $\Gamma_{n}$. Hence,
for all sufficiently large $n$, $\Gamma_{n}$ contains no non-trivial
element of translation length at most $2R$. Therefore $\mathrm{InjRad}\left(X_{n}\right)>R$.
This implies 
\[
\frac{\mathrm{Vol}\left\{ z\in X_{n}:\mathrm{InjRad}_{X_{n}}(z)<R\right\} }{\mathrm{Vol}\left(X_{n}\right)}\to0.
\]
Since $R$ is arbitrary, $X_{n}$ converges to $\mathbb{D}$ in the
Benjamini-Schramm sense. Let 
\[
p_{n}:X_{n}\to X_{0}
\]
be the covering map. Choose a non-zero real-valued Laplace eigenfunction
$a_{0}$ on $X_{0}$, normalized so that $\left\Vert a_{0}\right\Vert _{\infty}=1$,
with eigenvalue $\lambda=\frac{1+t^{2}}{4}>\frac{1}{4}$. Then we
set $a_{n}=a_{0}\circ p_{n}$, so that $a_{n}$ is also an eigenfunction
of the Laplacian with spectral parameter $t$. Then, (\ref{eq:I_e_2})
implies the identity contribution of $a_{n}$ is
\[
I_{e}=c_{k}(it)\left\Vert a_{n}\right\Vert _{2}^{2}.
\]
 Since $a_{n}=a_{0}\circ p_{n},$
\[
\left\Vert a_{n}\right\Vert _{2}^{2}=\left\Vert a_{0}\right\Vert _{2}^{2}\cdot\deg p_{n}=\left\Vert a_{0}\right\Vert _{2}^{2}\cdot\frac{\mathrm{Vol}\left(X_{n}\right)}{\mathrm{Vol}\left(X_{0}\right)}.
\]
Moreover, $\left(a_{n}\right)_{n\in\mathbb{N}}$ is also uniformly
bounded by 
\[
\left\Vert a_{n}\right\Vert _{\infty}=\left\Vert a_{0}\right\Vert _{\infty}=1.
\]
Then, by Corollary \ref{cor:BS nonidentity term} and (\ref{eq:HS=00003Di+n}),
the Hilbert-Schmidt norm of $Q_{k}A_{n}Q_{k}$, where $A_{n}$ is
the pull-back of $a_{n}$ to $\Gamma_{n}\backslash G$ as in (\ref{eq:A_n pull back}),
is 
\begin{equation}
\left\Vert Q_{k}A_{n}Q_{k}\right\Vert _{\mathrm{HS}}^{2}=\left(\frac{k-1}{4\pi}\right)^{2}c_{k}\left(it\right)\left\Vert a_{0}\right\Vert _{2}^{2}\cdot\frac{\mathrm{Vol}\left(X_{n}\right)}{\mathrm{Vol}\left(X_{0}\right)}+o_{k,r_{0}}\left(\mathrm{Vol}\left(X_{n}\right)\right).\label{eq:counterexample-HS}
\end{equation}
Since we assume $a_{0}$ to be real, $A_{n}$ is also real, and hence
$Q_{k}A_{n}Q_{k}$ is self-adjoint. Since $\sigma_{k}\left(\mathcal{M}_{k}\left(\Gamma_{n}\right)\right)$
is finite-dimensional, the spectral theorem implies $Q_{k}A_{n}Q_{k}$
is diagonalizable. We choose $\psi_{n,j}$ to be the eigenbasis of
$Q_{k}A_{n}Q_{k}$ so that 
\[
\sum_{j=1}^{d_{n,k}}\left|\left\langle \psi_{n,j},A_{n}\psi_{n,j}\right\rangle \right|^{2}=\left\Vert Q_{k}A_{n}Q_{k}\right\Vert _{\mathrm{HS}}^{2}.
\]
Since $a_{n}$ has a non-zero eigenvalue, $A_{n}$ is orthogonal to
the trivial representation. Equivalently, $A_{n}$ is mean zero. Then,
by (\ref{eq:d_n,k and Vol}) we have the quantum variance of $A_{n}$
with respect to $\left\{ \psi_{n,j}\right\} _{j=1}^{d_{n,k}}$ is
\[
\frac{1}{d_{n,k}}\sum_{j=1}^{d_{n,k}}\left|\left\langle \psi_{n,j},A_{n}\psi_{n,j}\right\rangle \right|^{2}=\frac{k-1}{4\pi}c_{k}(it)\frac{\left\Vert a_{0}\right\Vert _{2}^{2}}{\mathrm{Vol}\left(X_{0}\right)}+o_{k,r_{0}}\left(1\right).
\]
Hence, as $n\to\infty$,
\[
\frac{1}{d_{n,k}}\sum_{j=1}^{d_{n,k}}\left|\left\langle \psi_{n,j},A_{n}\psi_{n,j}\right\rangle \right|^{2}\to\frac{k-1}{4\pi}c_{k}(it)\frac{\left\Vert a_{0}\right\Vert _{2}^{2}}{\mathrm{Vol}\left(X_{0}\right)}.
\]
This proves Theorem \ref{thm:bad basis}. 
\begin{rem}
\label{rem:uniform spectral gap}The tower above may moreover be chosen
to have a uniform spectral gap. Indeed, we may take $\Gamma_{0}$
to be a cocompact arithmetic lattice and $\left(\Gamma_{n}\right)$
to be a nested tower of finite-index congruence subgroups with trivial
intersection. By Selberg's $3/16$ theorem \cite{key-17}, together
with the Jacquet-Langlands correspondence \cite{key-16}, the corresponding
compact hyperbolic surfaces $X_{n}=\Gamma_{n}\backslash\mathbb{D}$
satisfy 
\[
\lambda_{1}\left(X_{n}\right)\geq\frac{3}{16}
\]
uniformly in $n$. Thus the counterexample persists even under a uniform
spectral gap assumption.
\end{rem}

\section{Expectation Over Haar-Random Bases}\label{sec:Expectation-Over-Haar-Random}

For a given $\Gamma$, let 
\[
d=\mathrm{dim}\mathcal{M}_{k}\left(\Gamma\right).
\]
Fix an orthonormal basis of $\sigma_{k}\left(\mathcal{M}_{k}\left(\Gamma\right)\right)$
\[
\mathcal{E}=\left\{ e_{1},\ldots,e_{d}\right\} .
\]
For $U\in\mathrm{SU}(d)$, let 
\[
\psi_{j}^{U}=Ue_{j},\qquad1\leq j\leq d,
\]
and define
\[
\mathcal{V}_{\Gamma,k}\left(U,A\right)=\frac{1}{d}\sum_{j=1}^{d}\left|\left\langle \psi_{j}^{U},A\psi_{j}^{U}\right\rangle -\fint A\right|^{2}.
\]
Equip $\mathrm{SU}\left(d\right)$ with normalized Haar probability
measure.
\begin{prop}
\label{prop:general meanQV}Let $\mathcal{H}$ be a finite-dimensional
complex Hilbert space with $\mathrm{dim}\mathcal{H}=d$ and let $T:\mathcal{H}\to\mathcal{H}$
be a linear operator. Fix an orthonormal basis
\[
\mathcal{E}=\left\{ e_{1},\ldots e_{d}\right\} .
\]
For $U\in\mathrm{SU}(\mathcal{H}),$ by equipping $\mathrm{SU}\left(\mathcal{H}\right)$
with normalized Haar probability measure, we have 
\[
\mathbb{E}_{U}\left[\frac{1}{d}\sum_{j=1}^{d}\left|\left\langle Ue_{j},TUe_{j}\right\rangle \right|^{2}\right]=\frac{\left\Vert T\right\Vert _{\mathrm{HS}}^{2}+\left|\mathrm{Tr}T\right|^{2}}{d\left(d+1\right)}.
\]
\end{prop}
\begin{proof}
The case $d=1$ is immediate, so we assume $d\geq2$. For each fixed
$j$, the random vector $Ue_{j}$ is uniformly distributed on the
unit sphere of $\mathcal{H}$. Since both $\mathrm{U}\left(d\right)$
and $\mathrm{SU}\left(d\right)$ act transitively on the unit sphere
for $d\geq2$, a Haar-distributed matrix from either group sends any
fixed unit vector to a uniformly distributed unit vector. Let $v$
be a uniformly distributed unit vector and write 
\[
v=\sum_{a=1}^{d}v_{a}e_{a}.
\]
Then 
\[
\left\langle v,Tv\right\rangle =\sum_{a,b=1}^{d}T_{ab}\overline{v_{a}}v_{b}.
\]
 Hence 
\[
\left|\left\langle v,Tv\right\rangle \right|^{2}=\sum_{a,b,c,e}T_{ab}\overline{T_{ce}}v_{b}\overline{v_{a}}v_{c}\overline{v_{e}}.
\]
By the Weingarten formula of Collins and Śniady \cite[Corollary 2.4]{key-18},
the fourth moments of $v$ satisfy
\[
\mathbb{E}\left[v_{b}\overline{v_{a}}v_{c}\overline{v_{e}}\right]=\frac{\delta_{ab}\delta_{ce}+\delta_{be}\delta_{ca}}{d(d+1)}.
\]
It follows that 
\begin{align*}
\mathbb{E}_{v}\left[\left|\left\langle v,Tv\right\rangle \right|^{2}\right] & =\frac{1}{d(d+1)}\sum_{a,b,c,e}T_{ab}\overline{T_{ce}}\left(\delta_{ab}\delta_{ce}+\delta_{be}\delta_{ca}\right)\\
 & =\frac{1}{d(d+1)}\left[\sum_{a,c}T_{aa}\overline{T_{cc}}+\sum_{a,b}T_{ab}\overline{T_{ab}}\right]\\
 & =\frac{\left|\mathrm{Tr}T\right|^{2}+\left\Vert T\right\Vert _{\mathrm{HS}}^{2}}{d(d+1)}.
\end{align*}
Therefore, by linearity of expectation, we have 
\[
\mathbb{E}_{U}\left[\frac{1}{d}\sum_{j=1}^{d}\left|\left\langle Ue_{j},TUe_{j}\right\rangle \right|^{2}\right]=\frac{\left\Vert T\right\Vert _{\mathrm{HS}}^{2}+\left|\mathrm{Tr}T\right|^{2}}{d\left(d+1\right)}.
\]
\end{proof}
\begin{cor}
\label{cor:EVO(1/d)}Let 
\[
X_{n}=\Gamma_{n}\backslash\mathbb{D}
\]
be a sequence of surfaces converging to $\mathbb{D}$ in the Benjamini-Schramm
sense, and suppose 
\[
\mathrm{InjRad}\left(X_{n}\right)\geq r_{0}>0
\]
uniformly in $n$. Let 
\[
\rho_{n}=\inf_{R>0}\left[\sum_{m\geq2R}e^{\left(1-k/2\right)m}+\frac{\mathrm{Vol}\left(\left\{ z\in X_{n}:\mathrm{InjRad}_{X_{n}}z\leq R\right\} \right)}{\mathrm{Vol\left(X_{n}\right)}}\right].
\]
Then $\rho_{n}\to0$ as $n\to\infty.$ Moreover, uniformly for all
mean-zero observables $A_{n}$ satisfying
\[
\left\Vert A_{n}\right\Vert _{\infty}\leq1,
\]
we have
\[
\mathbb{E}_{U}\left[\mathcal{V}_{n,k}\left(U,A_{n}\right)\right]\lesssim_{k,r_{0}}\frac{1}{d_{n,k}}+\rho_{n}^{2}.
\]
Here $\mathcal{V}_{n,k}\left(U,A_{n}\right)=\mathcal{V}_{\Gamma_{n},k}\left(U,A_{n}\right).$
\end{cor}
\begin{proof}
We apply Proposition \ref{prop:general meanQV} to the case 
\[
\mathcal{H}=\sigma_{k}\left(\mathcal{M}_{k}\left(\Gamma_{n}\right)\right),\qquad T=Q_{k}A_{n}Q_{k}.
\]
Then 
\[
\mathbb{E}_{U}\left[\mathcal{V}_{n,k}\left(U,A_{n}\right)\right]=\frac{\left\Vert Q_{k}A_{n}Q_{k}\right\Vert _{\mathrm{HS}}^{2}+\left|\mathrm{Tr}\left(Q_{k}A_{n}Q_{k}\right)\right|^{2}}{d_{n,k}\left(d_{n,k}+1\right)}.
\]
We write the Hilbert-Schmidt norm term as in (\ref{eq:HS=00003Di+n}).
For the non-identity term, (\ref{eq:d_n,k and Vol}) and Corollary
\ref{cor:BS nonidentity term} give $o_{k,r_{0}}\left(d_{n,k}\right)$.
For the identity term, since $\sigma_{s}$ unitary and $\left\Vert e_{s}\right\Vert =1$,
\[
\left|\left\langle e_{s},\sigma_{s}(g)e_{s}\right\rangle \right|\leq1.
\]
Then by considering the formula for $c_{k}(s)$ in (\ref{eq:c_k_G}),
we have 
\begin{align*}
\left|c_{k}(s)\right| & \leq\int_{G}h_{k}(g)dg\\
 & =\int_{G}q_{k}(g)q_{k}\left(g^{-1}\right)dg.
\end{align*}
Since $Q_{k}^{2}=Q_{k}$, by comparing the kernel, we have 
\[
\frac{k-1}{4\pi}\int_{G}q_{k}\left(z^{-1}x\right)q_{k}\left(y^{-1}z\right)dz=q_{k}\left(y^{-1}x\right).
\]
Then, by taking $y=x=e$, we have 
\begin{equation}
\left|c_{k}(s)\right|\leq\frac{4\pi}{k-1}.\label{eq:c_k general bound}
\end{equation}
Then (\ref{eq:I_e_2}) implies that 
\[
I_{e}\leq\frac{4\pi}{k-1}\left\Vert A_{n}\right\Vert _{2}^{2}\leq\frac{4\pi}{k-1}\mathrm{Vol}\left(X_{n}\right).
\]
Therefore, (\ref{eq:d_n,k and Vol}) implies that 
\[
\left\Vert Q_{k}A_{n}Q_{k}\right\Vert _{\mathrm{HS}}^{2}=O_{k,r_{0}}\left(d_{n,k}\right).
\]
For the trace term, we can express it as 
\begin{align*}
\mathrm{Tr}\left(Q_{k}A_{n}Q_{k}\right) & =\int_{\Gamma_{n}\backslash G}A_{n}(x)B_{k}^{\Gamma_{n}}(x,x)dx\\
 & =\frac{k-1}{4\pi}\int_{\Gamma_{n}\backslash G}A_{n}\left(x\right)\sum_{\gamma\in\Gamma_{n}}q_{k}\left(x^{-1}\gamma x\right)dx.
\end{align*}
Since $A_{n}(x)$ is mean-zero, it is orthogonal to the trivial representation.
Since $Q_{k}$ acts as identity on $\sigma_{k}\left(\mathcal{M}_{k}\left(\Gamma\right)\right)$
and annihilates $\left(\sigma_{k}\left(\mathcal{M}_{k}\left(\Gamma\right)\right)\right)^{\perp}$,
we have 
\[
B_{k}^{\Gamma}(x,y)=\sum_{j=1}^{d_{k}}\psi_{j}(x)\overline{\psi_{j}(y)}.
\]
The component of $B_{k}^{\Gamma_{n}}(x,x)$ in trivial representation
is 
\[
\frac{1}{\mathrm{Vol}\left(X_{n}\right)}\int_{\Gamma_{n}\backslash G}B_{k}^{\Gamma_{n}}(x,x)dx=\frac{1}{\mathrm{Vol}\left(X_{n}\right)}\int_{\Gamma_{n}\backslash G}\sum_{j=1}^{d_{n,k}}\left|\psi_{n,j}\left(x\right)\right|^{2}dx=\frac{d_{n,k}}{\mathrm{Vol}\left(X_{n}\right)}=\frac{k-1}{4\pi}.
\]
Then we have
\begin{align*}
\mathrm{Tr}\left(Q_{k}A_{n}Q_{k}\right) & =\int_{\Gamma_{n}\backslash G}A_{n}(x)\left(\frac{k-1}{4\pi}\sum_{\gamma\in\Gamma_{n}}q_{k}\left(x^{-1}\gamma x\right)-\frac{k-1}{4\pi}\right)dx\\
 & =\frac{k-1}{4\pi}\int_{\Gamma_{n}\backslash G}A_{n}(x)\sum_{\gamma\not=e}q_{k}\left(x^{-1}\gamma x\right)dx.
\end{align*}
Using the upper bound for $q_{k}$ in (\ref{eq:q_k upper bound}),
\[
\left|\mathrm{Tr}\left(Q_{k}A_{n}Q_{k}\right)\right|\lesssim_{k}\left\Vert A_{n}\right\Vert _{\infty}\int_{F_{n}}\sum_{\gamma\not=e}e^{-kr_{\gamma}(z)/2}d\mathrm{Vol}_{\mathbb{D}}(z).
\]
 Let $D(R)$ be the same as in Proposition \ref{prop:bound-nonidentity}.
Then for $z\in D(R)$, 
\begin{equation}
\sum_{\gamma\not=e}e^{-kr_{\gamma}(z)/2}\lesssim_{k,r_{0}}\sum_{m\geq2R}e^{\left(1-k/2\right)m}.\label{eq:r_=00005Cgamma}
\end{equation}
This is convergent as $k>2$. Then 
\[
\int_{D(R)}\sum_{\gamma\not=e}e^{-kr_{\gamma}(z)/2}d\mathrm{Vol}_{\mathbb{D}}(z)\lesssim_{k,r_{0}}\mathrm{Vol}\left(X_{n}\right)\sum_{m\geq2R}e^{\left(1-k/2\right)m}.
\]
For $z\in D(R)^{c}$, we have the general bound
\begin{equation}
\sum_{\gamma\not=e}e^{-kr_{\gamma}(z)/2}\lesssim_{k,r_{0}}\sum_{m\geq0}e^{\left(1-k/2\right)m}.\label{eq:Generalbound}
\end{equation}
Therefore 
\[
\int_{D(R)^{c}}\sum_{\gamma\not=e}e^{-kr_{\gamma}(z)/2}d\mathrm{Vol}_{\mathbb{D}}(z)\lesssim_{k,r_{0}}\mathrm{Vol}\left(D(R)^{c}\right).
\]
Therefore, 
\[
\left|\mathrm{Tr}\left(Q_{k}A_{n}Q_{k}\right)\right|\lesssim_{k,r_{0}}\left\Vert A_{n}\right\Vert _{\infty}\left[\mathrm{Vol}\left(X_{n}\right)\sum_{m\geq2R}e^{\left(1-k/2\right)m}+\mathrm{Vol}\left(\left\{ z\in X_{n}:\mathrm{InjRad}_{X_{n}}(z)\leq R\right\} \right)\right].
\]
Taking the infimum over $R>0$, we obtain 
\[
\left|\mathrm{Tr}\left(Q_{k}A_{n}Q_{k}\right)\right|\lesssim_{k,r_{0}}\mathrm{\mathrm{Vol}}\left(X_{n}\right)\rho_{n}\lesssim_{k,r_{0}}d_{n,k}\rho_{n}.
\]
Therefore 
\[
\mathbb{E}_{U}\left[\mathcal{V}_{n,k}\left(U,A_{n}\right)\right]\lesssim_{k,r_{0}}\frac{1}{d_{n,k}}+\rho_{n}^{2},
\]
which proves the result.
\end{proof}

\section{Concentration Over Haar-Random Bases}\label{sec:Concentration-Over-Haar-Random}

Corollary \ref{cor:EVO(1/d)} implies that for every sequence of observables
$\left(A_{n}\right)_{n\in\mathbb{N}}$, there exists a sequence of
orthonormal bases such that the corresponding quantum variance converges
to $0$. It is therefore natural to ask whether there exists a single
sequence of orthonormal bases for which the quantum variance converges
to $0$ simultaneously for every sequence of observables $\left(A_{n}\right)_{n\in\mathbb{N}}$.
In this section, we will show that this is indeed the case. We first
establish a general estimate that depends on the choice of parameters.
By choosing suitable parameters, we then prove Theorem \ref{thm:most coverge}.
\begin{prop}
\label{prop:main thm}Let $X_{n}=\Gamma_{n}\backslash\mathbb{D}$
be a sequence of compact hyperbolic surfaces converging to $\mathbb{D}$
in the Benjamini-Schramm sense and 
\[
\mathrm{InjRad}\left(X_{n}\right)\geq r_{0}>0
\]
uniformly in $n$. Fix $k>2$. There exist constants $C_{1},C_{2}>0$,
depending only on $k$ and $r_{0}$, such that for every choice of
\[
T_{n}>0,\qquad\delta_{n}>0,\qquad\varepsilon_{n}>0,
\]
there exists a measurable set $\Omega_{n}\subset\mathrm{SU}\left(d_{n,k}\right)$
satisfying 
\[
\mu_{n,k}\left(\Omega_{n}\right)\geq1-\sqrt{2}\left(1+\frac{2}{\delta_{n}}\right)^{2m_{n}(T_{n})}\exp\left(-\frac{d_{n,k}\varepsilon_{n}^{2}}{4}\right),
\]
where $m_{n}(T_{n})$ denotes the number of the non-zero eigenvalues
on $X_{n}$, counted with multiplicity, satisfying $\lambda\leq1/4+T_{n}^{2}$.
Moreover, there exists a sequence $\left(\tau_{n}\right)$ with $\tau_{n}\to0$
as $n\to\infty$ such that for every $U\in\Omega_{n}$ and every mean-zero
observable $A_{n}$ satisfying $\left\Vert A_{n}\right\Vert _{\infty}\leq1$,
\[
\mathcal{V}_{n,k}\left(U,A_{n}\right)\leq2C_{1}\left(\frac{1}{d_{n,k}}+\rho_{n}^{2}\right)+\frac{8\varepsilon_{n}}{\sqrt{d_{n,k}}}+2C_{2}\left(\delta_{n}^{2}+\left(1+T_{n}^{2k-2}\right)e^{-\pi T_{n}/2}+\tau_{n}\right).
\]
\end{prop}
Proposition \ref{prop:main thm} gives a quantitative estimate for
the quantum variance in terms of three auxiliary parameters $T_{n}$,
$\delta_{n}$, $\varepsilon_{n}$. To deduce Theorem \ref{thm:most coverge},
we must choose these parameters so that the Haar measure of $\Omega_{n}$
tends to one, while the upper bound for the quantum variance tends
to zero uniformly over all bounded mean-zero observables. Since $\rho_{n},\tau_{n}\to0$,
it suffices that
\[
2m_{n}\left(T_{n}\right)\log\left(1+\frac{2}{\delta_{n}}\right)-\frac{d_{n,k}\varepsilon_{n}^{2}}{4}\to-\infty,
\]
and 
\[
\delta_{n}\to0,\qquad T_{n}\to\infty,\qquad\frac{\varepsilon_{n}}{\sqrt{d_{n,k}}}\to0.
\]
To choose suitable parameters, we need to control $m_{n}\left(T\right)$
in terms of $d_{n,k}$. We adapt the proof of Weyl's law in Bergeron
\cite[§5.4.1]{key-15} by modifying the closed geodesic counting
estimate in Lemma 5.22. This gives the following bound of $m_{n}\left(T\right)$.
\begin{lem}
\label{lem:uniform bound for mnT}Let $X_{n}=\Gamma_{n}\backslash\mathbb{D}$
be a sequence of compact hyperbolic surfaces satisfying 
\[
\mathrm{InjRad}\left(X_{n}\right)\geq r_{0}>0
\]
 for all $n$ in $\mathbb{N}$. Then 
\[
m_{n}\left(T\right)\lesssim_{r_{0}}\mathrm{Vol}\left(X_{n}\right)\left(1+T^{2}\right),
\]
as $T\to\infty$, where $m_{n}\left(T\right)$ is defined as above.
\end{lem}
\begin{proof}
We use the packing argument from Proposition \ref{prop:bound-nonidentity}.
Since $\mathrm{InjRad}\left(X_{n}\right)\geq r_{0}>0$, we have for
each $z\in\mathbb{D}$ and $n\in\mathbb{N}$,
\begin{equation}
\#\left\{ \gamma\in\Gamma_{n}:d_{\mathbb{D}}\left(z,\gamma z\right)\leq L\right\} =O_{r_{0}}\left(e^{L}\right).\label{eq:count geodesic at fixed point}
\end{equation}
We use this estimate to obtain a uniform version of Bergeron's Lemma
5.22. Choose a maximal $r_{0}/2$-separated set of points $\left\{ p_{1},\ldots,p_{N_{n}}\right\} $in
$X_{n}$. Since the balls of radius $r_{0}/4$ centered at these points
are disjoint, we have 
\[
N_{n}\lesssim_{r_{0}}\mathrm{Vol}\left(X_{n}\right).
\]
Fix a lift $\tilde{p_{i}}\in\mathbb{D}$ of each $p_{i}$. Every oriented
closed geodesic $\gamma$ of length at most $L$ passes within distance
$r_{0}/2$ of some $p_{i}$. Choose such a $p_{i}$ and a shortest
path joining $p_{i}$ to $\gamma$. Following this path, traversing
$\gamma$, and returning along the reverse path gives a loop at $p_{i}$,
freely homotopic to $\gamma$, of length at most $L+r_{0}$. Its lift
starting at $\tilde{p_{i}}$ ends at $g_{\gamma}\tilde{p_{i}}$ for
some $g_{\gamma}\in\Gamma_{n}$, so $d_{\mathbb{D}}\left(\tilde{p_{i}},g_{\gamma}\tilde{p_{i}}\right)\leq L+r_{0}$.
For each fixed $i$, geodesics in distinct free homotopy classes give
elements $g_{\gamma}$ in distinct conjugacy classes, and these elements
are therefore distinct. Applying (\ref{eq:count geodesic at fixed point})
at each $\tilde{p_{i}}$ and summing over $i$, we obtain
\begin{equation}
\#\left\{ \gamma\in G\left(X_{n}\right):\ell(\gamma)\leq L\right\} \lesssim_{r_{0}}\mathrm{Vol}\left(X_{n}\right)e^{L},\label{eq:universal bound of geodesic}
\end{equation}
where $G\left(X_{n}\right)$ denotes the oriented closed geodesics,
including iterates.

We now follow the argument in Bergeron. For $\delta>0$, we apply
the Selberg trace formula to admissible pair 
\[
h(r)=e^{-\delta r^{2}},\qquad g(x)=\left(4\pi\delta\right)^{-1/2}e^{-x^{2}/4\delta},
\]
to obtain 
\begin{equation}
\sum_{j=0}^{\infty}e^{-\delta r_{n,j}^{2}}=\frac{\mathrm{Vol}\left(X_{n}\right)}{4\pi}\int_{-\infty}^{+\infty}re^{-\delta r^{2}}\tanh\left(\pi r\right)dr+\left(4\pi\delta\right)^{-1/2}\sum_{\gamma\in G\left(X_{n}\right)}\frac{\ell\left(\gamma_{0}\right)}{2\sinh\left(\ell\left(\gamma\right)/2\right)}e^{-\ell(\gamma)^{2}/4\delta},\label{eq:selberg trace formula}
\end{equation}
where $\lambda_{n,j}=1/4+r_{n,j}^{2}$, and $\gamma_{0}$ is the primitive
closed geodesic underlying $\gamma$. Since $\ell(\gamma)\geq2r_{0}$,
the uniform geodesic counting estimate in (\ref{eq:universal bound of geodesic})
gives 
\[
\left(4\pi\delta\right)^{-1/2}\sum_{\gamma\in G\left(X_{n}\right)}\frac{\ell\left(\gamma_{0}\right)}{2\sinh\left(\ell\left(\gamma\right)/2\right)}e^{-\ell(\gamma)^{2}/4\delta}=O_{r_{0}}\left(\mathrm{Vol}\left(X_{n}\right)\delta^{-1/2}e^{-r_{0}^{2}/2\delta}\right),
\]
which converges to $0$ as $\delta\to0$ for each fixed $n$. As in
Bergeron, integration by parts gives 
\[
\int_{-\infty}^{+\infty}re^{-\delta r^{2}}\tanh\left(\pi r\right)dr=\frac{1}{\delta}+O\left(1\right).
\]
Hence, by summing up these two terms, we can rewrite (\ref{eq:selberg trace formula})
as 
\[
\frac{1}{\mathrm{Vol}\left(X_{n}\right)}\sum_{j=0}^{\infty}e^{-\delta r_{n,j}^{2}}=\frac{1}{4\pi\delta}+O_{r_{0}}\left(1\right).
\]
Since we do not need an asymptotic formula for $m_{n}\left(T\right)$
as in Weyl's law, instead of applying the Tauberian Theorem of Karamata,
we use the following argument
\[
m_{n}\left(T\right)e^{-\delta T^{2}}\leq\sum_{j=0}^{\infty}e^{-\delta r_{n,j}^{2}},
\]
 then, for sufficiently small $\delta>0$,
\[
\frac{m_{n}\left(T\right)}{\mathrm{Vol}\left(X_{n}\right)}\leq e^{\delta T^{2}}\left(\frac{1}{4\pi\delta}+O_{r_{0}}\left(1\right)\right)
\]
 We now choose $\delta=\frac{1}{1+T^{2}}$, so that $\delta\to0$
as $T\to\infty$ and obtain that 
\[
\frac{m_{n}\left(T\right)}{\mathrm{Vol}\left(X_{n}\right)}\leq e^{\frac{T^{2}}{1+T^{2}}}\left(\frac{1+T^{2}}{4\pi}+O_{r_{0}}\left(1\right)\right)\leq e\left(\frac{1+T^{2}}{4\pi}+O_{r_{0}}\left(1\right)\right)\lesssim_{r_{0}}1+T^{2}.
\]
This proves the desired bound.
\end{proof}
We can now prove Theorem \ref{thm:most coverge} by choosing appropriate
parameters.
\begin{proof}[Proof of Theorem 1.2]
By Lemma \ref{lem:uniform bound for mnT} and the dimension formula
in (\ref{eq:d_n,k and Vol}), we obtain
\[
m_{n}\left(T\right)\lesssim_{k,r_{0}}d_{n,k}\left(1+T^{2}\right)
\]
uniformly in $n$. One possible choice of parameters is 
\[
T_{n}=d_{n,k}^{1/8},\qquad\delta_{n}=d_{n,k}^{-1/8},\qquad\varepsilon_{n}=d_{n,k}^{1/4}.
\]
Indeed in this case,
\begin{equation}
2m_{n}\left(T_{n}\right)\log\left(1+\frac{2}{\delta_{n}}\right)\lesssim_{k,r_{0}}d_{n,k}^{5/4}\log d_{n,k},\label{eq:a bound}
\end{equation}
whereas 
\[
\frac{d_{n,k}\varepsilon_{n}^{2}}{4}=\frac{1}{4}d_{n,k}^{3/2}.
\]
Hence, for all sufficiently large $n$,
\begin{equation}
2m_{n}\left(T_{n}\right)\log\left(1+\frac{2}{\delta_{n}}\right)-\frac{d_{n,k}\varepsilon_{n}^{2}}{4}\leq C_{k,r_{0}}d_{n,k}^{5/4}\log d_{n,k}-\frac{1}{4}d_{n,k}^{3/2}\leq-\frac{1}{8}d_{n,k}^{3/2},\label{eq:exponential upper bound}
\end{equation}
where here $C_{k,r_{0}}$ is the constant depending on $k$ and $r_{0}$
in (\ref{eq:a bound}). The latter goes to $-\infty$ as $n\to\infty$.
Hence 
\[
\mu_{n,k}\left(\Omega_{n}\right)\to1.
\]
Moreover, 
\[
\frac{\varepsilon_{n}}{\sqrt{d_{n,k}}}=d_{n,k}^{-1/4},\qquad\delta_{n}^{2}=d_{n,k}^{-1/4},
\]
while 
\[
\left(1+T_{n}^{2k-2}\right)e^{-\pi T_{n}/2}\to0.
\]
Since $\rho_{n}\to0$ and $\tau_{n}\to0$, the upper bound in Proposition
\ref{prop:main thm} tends to zero. For an arbitrary observable $a_{n}\in L^{\infty}\left(X_{n}\right)$
with $\left\Vert a_{n}\right\Vert _{\infty}\leq1$, set 
\[
b_{n}=\frac{1}{2}\left(a_{n}-\fint_{X_{n}}a_{n}d\mathrm{Vol}_{X_{n}}\right).
\]
Then $b_{n}$ has mean zero and $\left\Vert b_{n}\right\Vert _{\infty}\leq1$.
Writing $B_{n}$ for its lift to $\Gamma_{n}\backslash G$, we have
\[
\mathcal{V}_{n,k}\left(U,a_{n}\right)=4\mathcal{V}_{n,k}\left(U,B_{n}\right).
\]
Applying Proposition \ref{prop:main thm} to $B_{n}$, we therefore
obtain 
\[
\sup_{\left\Vert a_{n}\right\Vert _{\infty}\leq1}\mathcal{V}_{n,k}\left(U,a_{n}\right)\to0
\]
for all $U\in\Omega_{n}$. Together with $\mu_{n,k}\left(\Omega_{n}\right)\to1$,
this proves Theorem \ref{thm:most coverge}.
\end{proof}
\begin{cor}[Almost sure convergence]
Under the assumptions of Theorem \ref{thm:most coverge}, suppose
in addition that 
\[
\sum_{n=1}^{\infty}\exp\left(-\frac{1}{8}d_{n,k}^{3/2}\right)<\infty,
\]
Let $\left(U_{n}\right)_{n\in\mathbb{N}}$ be a sequence of independent
Haar-distributed random matrices, with $U_{n}\in\mathrm{SU}\left(d_{n,k}\right)$.
Then
\[
\sup_{\left\Vert A_{n}\right\Vert _{\infty}\leq1}\mathcal{V}_{n,k}\left(U_{n},A_{n}\right)\to0
\]
 almost surely with respect to the product Haar measure, where the
supremum is taken over all observables $A_{n}$ obtained as lifts
of $a_{n}\in L^{\infty}\left(X_{n}\right)$, as defined in (\ref{eq:A_n pull back}).
\end{cor}
This follows directly from the choice of parameter above and the first
Borel-Cantelli lemma.

It remains to prove Proposition \ref{prop:main thm}. We begin by
decomposing observables into low- and high-frequency components, and
then establish the required concentration estimates for Haar-random
bases.

\subsection{Low- and high-frequency decomposition }

Let $X=\Gamma\backslash\mathbb{D}$ be a compact hyperbolic surface
satisfying $\mathrm{InjRad}\left(X\right)\geq r_{0}>0$, and let $A\in L^{\infty}\left(\Gamma\backslash G\right)^{K}$
be a mean-zero observable. Using the spherical decomposition in (\ref{eq:decomp of A}),
we may write 
\[
A=A_{\mathrm{comp}}+\sum_{j}A_{it_{j}},
\]
where $A_{\mathrm{comp}}$ denotes the sum of the complementary-series
components and $A_{it_{j}}$ denotes the component corresponding to
the spherical principal-series representation with spectral parameter
$t_{j}.$ 

For $T>0$, let $\mathscr{P}(T)$ denote the subspace of $L^{2}\left(\Gamma\backslash G\right)^{K}$
consisting of all complementary-series components together with principal-series
components with spectral parameter $\left|t\right|\leq T$. Let 
\[
P_{T}:L^{2}\left(\Gamma\backslash G\right)^{K}\to\mathscr{P}(T)
\]
 be the $L^{2}$-orthogonal projection. Then,
\[
P_{T}A=A_{\mathrm{comp}}+\sum_{\left|t_{j}\right|\leq T}A_{it_{j}},\qquad\left(1-P_{T}\right)A=\sum_{\left|t_{j}\right|>T}A_{it_{j}}.
\]
Then 
\begin{equation}
A=P_{T}A+\left(1-P_{T}\right)A.\label{eq:T decomp}
\end{equation}
The Hilbert-Schmidt norm calculated in Section \ref{sec:The-Hilbert-Schmidt-Norm}
gives us 
\[
\mathcal{V}_{\Gamma,k}\left(U,A\right)\lesssim_{k}\frac{1}{d_{k}}\left(I_{e}+\sum_{\gamma\not=e}I_{\gamma}\right).
\]
By Proposition \ref{prop:bound-nonidentity}, the non-identity contribution
satisfies
\[
\left|\sum_{\gamma\not=e}I_{\gamma}\right|\lesssim_{k,r_{0}}\left\Vert A\right\Vert _{\infty}^{2}\left[\sum_{m\geq2R}(m+1)e^{m-\frac{k}{2}m}\cdot\mathrm{Vol}(X)+\mathrm{Vol}\left(\left\{ z\in F:\mathrm{InjRad}_{X}(z)\leq R\right\} \right)\right]
\]
By (\ref{eq:I_e_2}), the identity contribution is
\[
I_{e}=\sum_{\sigma_{s}\not\simeq\mathbf{1}}\sum_{a=1}^{m\left(\sigma_{s},\Gamma\right)}c_{k}(s)\left\Vert A_{\sigma_{s},a}\right\Vert ^{2}.
\]
For every spherical unitary representation occurring in the low-frequency
space, by (\ref{eq:c_k general bound}),
\[
\left|c_{k}(s)\right|\leq\frac{4\pi}{k-1}.
\]
We then calculate an upper bound for $c_{k}(it)$ for $t>T$. Since
from Proposition \ref{prop:c_k calculation }, we have 
\[
c_{k}(it)=4\pi\frac{\left|\Gamma\left(k-\frac{1}{2}+\frac{it}{2}\right)\right|^{2}}{\Gamma(k)^{2}}.
\]
Using the recurrence relation for Gamma function and 
\[
\left|\Gamma\left(\frac{1}{2}+\frac{it}{2}\right)\right|^{2}=\frac{\pi}{\cosh\left(\pi t/2\right)},
\]
we obtain 
\[
c_{k}(it)=\frac{4^{2-k}\pi^{2}}{\left(\left(k-1\right)!\right)^{2}}\frac{\prod_{j=0}^{k-2}\left(t^{2}+\left(2j+1\right)^{2}\right)}{\cosh\left(\pi t/2\right)}.
\]
Hence, for $\left|t\right|>T$, there exists $C_{k}$ such that
\[
\left|c_{k}(it)\right|\leq C_{k}\left(1+T^{2k-2}\right)e^{-\pi T/2}=\eta_{k}\left(T\right).
\]
In particular 
\[
\eta_{k}\left(T\right)\to0\qquad\text{as }T\to\infty.
\]
Using the decomposition in (\ref{eq:T decomp}), we can bound the
identity term by 
\[
I_{e}\lesssim_{k}\left\Vert P_{T}A\right\Vert _{2}^{2}+\eta_{k}\left(T\right)\left\Vert A\right\Vert _{2}^{2}.
\]
Thus by (\ref{eq:d_n,k and Vol}), 
\begin{equation}
\mathcal{V}_{\Gamma,k}\left(U,A\right)\lesssim_{k,r_{0}}\frac{\left\Vert P_{T}A\right\Vert _{2}^{2}}{\mathrm{Vol}(X)}+\eta_{k}(T)\left\Vert A\right\Vert _{\infty}^{2}+\mathrm{GeomErr}_{X}\left(R\right)\left\Vert A\right\Vert _{\infty}^{2}.\label{eq:HS upperbound T decomp}
\end{equation}
Here
\[
\mathrm{GeomErr}_{X}\left(R\right)=\frac{1}{\mathrm{Vol}(X)}\left(\sum_{m\geq2R}(m+1)e^{m-\frac{k}{2}m}\cdot\mathrm{Vol}(X)+\mathrm{Vol}\left(\left\{ z\in F:\mathrm{InjRad}_{X}(z)\leq R\right\} \right)\right)
\]
Thus, by Benjamini-Schramm convergence, 
\[
\inf_{R>0}\mathrm{GeomErr}_{X_{n}}\left(R\right)\to0,\qquad\text{as }n\to\infty.
\]
Moreover, $\eta_{k}(T)\to0$ as $T\to\infty$. Hence, up to these
two error terms, the quantum variance is controlled by the $L^{2}$-norm
of the low-frequency component $P_{T}A$.
\begin{lem}
\label{lem:lowcontribution inhighT and continuity}Let $A,B\in L^{\infty}\left(\Gamma\backslash G\right)^{K}$
be mean-zero. Then 
\[
\left|\mathcal{V}_{\Gamma,k}\left(U,A\right)^{1/2}-\mathcal{V}_{\Gamma,k}\left(U,B\right)^{1/2}\right|\leq\mathcal{V}_{\Gamma,k}\left(U,A-B\right)^{1/2}.
\]
Consequently, if 
\[
\left\Vert A\right\Vert _{\infty},\left\Vert B\right\Vert _{\infty}\leq1,
\]
 then by applying the estimate (\ref{eq:HS upperbound T decomp})
to $A-B$,
\[
\left|\mathcal{V}_{\Gamma,k}\left(U,A\right)^{1/2}-\mathcal{V}_{\Gamma,k}\left(U,B\right)^{1/2}\right|\lesssim_{k,r_{0}}\left(\frac{\left\Vert P_{T}\left(A-B\right)\right\Vert _{2}^{2}}{\mathrm{Vol}(X)}+\eta_{k}(T)+\mathrm{GeomErr}_{X}\left(R\right)\right)^{1/2}.
\]
\end{lem}
\begin{proof}
The first inequality follows immediately from the triangle inequality
in $l^{2}$. Then using 
\[
\left\Vert A-B\right\Vert _{\infty}\leq2
\]
 gives the result.
\end{proof}

\subsection{Gromov-Milman concentration}

We next obtain a concentration estimate for quantum variance associated
with a fixed bounded observable. This can be achieved using the concentration
theorem of Gromov and Milman \cite[§1.2, Remark 2]{key-3} together
with a finite net argument. 
\begin{thm}[Gromov-Milman]
Let $X$ be a closed connected Riemannian manifold, and let $\mu$
denote the normalized Riemannian volume measure on $X$. Define 
\[
R(X)\coloneqq\inf_{\left\Vert \tau\right\Vert =1}\mathrm{Ric}\left(\tau,\tau\right),
\]
where the infimum is taken over all unit tangent vectors $\tau$ on
$X$. Let $A\subset X$ be a measurable subset satisfying 
\[
\mu\left(A\right)\geq\frac{1}{2},
\]
 and denote by $N_{\varepsilon}(A)$ the $\varepsilon$-neighborhood
of $A$. If $R(X)>0$, then 
\[
\mu\left(N_{\varepsilon}(A)\right)\geq1-\sqrt{2}\exp\left(-\frac{\varepsilon^{2}R(X)}{2}\right).
\]
\end{thm}
We apply this theorem to $X=\mathrm{SU}\left(d_{n,k}\right)$, equipped
with the bi-invariant Riemannian metric induced by the Hilbert-Schmidt
inner product
\[
\left\langle X,Y\right\rangle =-\mathrm{Tr}\left(XY\right).
\]
The normalized Riemannian volume measure coincides with the normalized
Haar probability measure $\mu_{n,k}$. We set
\[
A=E_{A_{n}}=\left\{ U\in\mathrm{SU}\left(d_{n,k}\right):\mathcal{V}_{n,k}\left(U,A_{n}\right)\leq m_{A_{n}}\right\} ,
\]
where $m_{A_{n}}$ is the median of the function $U\mapsto\mathcal{V}_{n,k}\left(U,A_{n}\right)$.
In this case, 
\[
R\left(\mathrm{SU}\left(d_{n,k}\right)\right)=\frac{d_{n,k}}{2}.
\]
Then we have 
\begin{equation}
\mu_{n,k}\left(N_{\varepsilon}\left(E_{A_{n}}\right)\right)\geq1-\sqrt{2}\exp\left(-\frac{d_{n,k}\varepsilon^{2}}{4}\right).\label{eq:Gromov-Milman}
\end{equation}

\subsection{Lipschitz continuity of quantum variance in the choice of basis}
\begin{lem}
\label{lem:Lipinbasis}Let $U,V\in\mathrm{SU}\left(d_{k}\right)$,
and let $A\in L^{\infty}\left(\Gamma\backslash G\right)$ be a mean-zero
observable satisfying $\left\Vert A\right\Vert _{\infty}\leq1$,
\[
\left|\mathcal{V}_{\Gamma,k}\left(U,A\right)-\mathcal{V}_{\Gamma,k}\left(V,A\right)\right|\leq\frac{4}{\sqrt{d_{k}}}d_{\mathrm{SU}\left(d_{k}\right)}\left(U,V\right),
\]
where $d_{\mathrm{SU}\left(d_{k}\right)}$ is the Riemannian distance
associated with the Hilbert-Schmidt metric.
\end{lem}
\begin{proof}
For convenience, we write 
\[
u_{j}=\psi_{j}^{U},\qquad v_{j}=\psi_{j}^{V}.
\]
 Then 
\begin{align*}
\left|\mathcal{V}_{\Gamma,k}\left(U,A\right)-\mathcal{V}_{\Gamma,k}\left(V,A\right)\right| & \leq\frac{1}{d_{k}}\sum_{j=1}^{d_{k}}\left|\left|\left\langle u_{j},Au_{j}\right\rangle \right|^{2}-\left|\left\langle v_{j},Av_{j}\right\rangle \right|^{2}\right|\\
 & \leq\frac{2}{d_{k}}\sum_{j=1}^{d_{k}}\left|\left\langle u_{j},Au_{j}\right\rangle -\left\langle v_{j},Av_{j}\right\rangle \right|,
\end{align*}
where the last inequality results from the triangle inequality and
the fact that $x^{2}$ is Lipschitz on $\left[0,1\right]$ with the
Lipschitz constant to be $2$.

As 
\begin{align*}
\left|\left\langle u_{j},Au_{j}\right\rangle -\left\langle v_{j},Av_{j}\right\rangle \right| & \leq\left|\left\langle \left(u_{j}-v_{j}\right),Au_{j}\right\rangle \right|+\left|\left\langle v_{j},A\left(u_{j}-v_{j}\right)\right\rangle \right|\\
 & \leq\left\Vert A\right\Vert _{\infty}\left\Vert u_{j}\right\Vert \left\Vert u_{j}-v_{j}\right\Vert +\left\Vert A\right\Vert _{\infty}\left\Vert v_{j}\right\Vert \left\Vert u_{j}-v_{j}\right\Vert \\
 & \leq2\left\Vert u_{j}-v_{j}\right\Vert ,
\end{align*}
 we have 
\begin{align*}
\left|\mathcal{V}_{\Gamma,k}\left(U,A\right)-\mathcal{V}_{\Gamma,k}\left(V,A\right)\right| & \leq\frac{4}{d_{k}}\sum_{j=1}^{d_{k}}\left\Vert u_{j}-v_{j}\right\Vert \\
 & =\frac{4}{d_{k}}\sum_{j=1}^{d_{k}}\left\Vert \left(U-V\right)e_{j}\right\Vert \\
 & \leq\frac{4}{\sqrt{d_{k}}}\left(\sum_{j=1}^{d_{k}}\left\Vert \left(U-V\right)e_{j}\right\Vert ^{2}\right)^{1/2}\\
 & =\frac{4}{\sqrt{d_{k}}}\left\Vert U-V\right\Vert _{\mathrm{HS}}\\
 & \leq\frac{4}{\sqrt{d_{k}}}d_{\mathrm{SU}\left(d_{k}\right)}\left(U,V\right).
\end{align*}
\end{proof}
For $U\in N_{\varepsilon_{n}}\left(E_{A_{n}}\right)$, there exists
$V\in E_{A_{n}}$ such that 
\[
d_{\mathrm{SU}\left(d_{n,k}\right)}\left(U,V\right)<\varepsilon_{n}.
\]
Lemma \ref{lem:Lipinbasis} implies that 
\[
\mathcal{V}_{n,k}\left(U,A_{n}\right)\leq\mathcal{V}_{n,k}\left(V,A_{n}\right)+\frac{4}{\sqrt{d_{n,k}}}\varepsilon_{n}\leq m_{A_{n}}+\frac{4}{\sqrt{d_{n,k}}}\varepsilon_{n}.
\]
Therefore,
\begin{equation}
\left\{ U\in\mathrm{SU}\left(d_{n,k}\right):\mathcal{V}_{n,k}\left(U,A_{n}\right)\leq m_{A_{n}}+\frac{4}{\sqrt{d_{n,k}}}\varepsilon_{n}\right\} \supseteqq N_{\varepsilon_{n}}\left(E_{A_{n}}\right).\label{eq:Contain1}
\end{equation}

\subsection{A $\delta$-net for low-frequency observables}

Since $X$ is compact, $\mathscr{P}(T)$ is finite dimensional. Let
$m\left(T\right)$ denote the number of non-zero Laplace eigenvalues,
counted with multiplicity, satisfying $\lambda\leq\frac{1}{4}+T^{2}$.
We have 
\[
\mathrm{dim}_{\mathbb{C}}\mathscr{P}\left(T\right)=\#\left\{ j:0<\lambda_{j}\leq\frac{1}{4}+\frac{T^{2}}{4}\right\} \leq m\left(T\right).
\]
We consider the set 
\[
\mathscr{A}(T)=\left\{ P_{T}A:A\text{ is mean-zero and}\left\Vert A\right\Vert _{\infty}\leq1\right\} \subset\mathscr{P}(T).
\]
The following elementary covering lemma will allow us to reduce the
uniform estimate over bounded observables to finitely many observables.
\begin{lem}
\label{lem:net}For every $\delta>0$, there exists a finite set 
\[
\mathcal{N}\left(T,\delta\right)=\left\{ B_{1},\ldots,B_{N}\right\} \subset\mathscr{A}(T)
\]
such that, for every mean-zero $A$ with $\left\Vert A\right\Vert _{\infty}\leq1$,
there exists $B_{i}\in\mathcal{N}(T,\delta)$ satisfying 
\[
\left\Vert P_{T}A-B_{i}\right\Vert _{2}\leq\delta\mathrm{Vol}\left(X\right)^{1/2},\qquad X=\Gamma\backslash\mathbb{D}.
\]
Moreover, 
\[
N\leq\left(1+\frac{2}{\delta}\right)^{2m(T)}.
\]
\end{lem}
\begin{proof}
Since $P_{T}$ is an orthogonal projection on $L^{2}\left(\Gamma\backslash G\right)$,
\[
\left\Vert P_{T}A\right\Vert _{2}\leq\left\Vert A\right\Vert _{2}.
\]
For $\left\Vert A\right\Vert _{\infty}\leq1$,
\[
\left\Vert A\right\Vert _{2}^{2}\leq\mathrm{Vol}\left(X\right),
\]
and hence 
\[
\left\Vert P_{T}A\right\Vert _{2}\leq\mathrm{Vol}(X)^{1/2}.
\]
Therefore
\[
\mathscr{A}(T)\subset B_{\mathscr{P}(T)}\left(0,\mathrm{Vol}\left(X\right)^{1/2}\right).
\]
Choose a maximal $\delta\mathrm{Vol}(X)^{1/2}$-separated subset
\[
\mathcal{N}\left(T,\delta\right)=\left\{ B_{1},\ldots,B_{N}\right\} \subset\mathscr{A}(T).
\]
By maximality, $\mathcal{N}(T,\delta)$ is a $\delta\mathrm{Vol}(X)^{1/2}$-net
for $\mathscr{A}(T)$. Indeed, if there exists $B\in\mathscr{A}(T)$
such that 
\[
\left\Vert B-B_{i}\right\Vert _{2}>\delta\mathrm{Vol}(X)^{1/2}
\]
for every $i$, then $B$ could be added to $\mathcal{N}(T,\delta)$,
contradicting maximality. 

It remains to bound $N$. Since the points $B_{i}$ are $\delta\mathrm{Vol}(X)^{1/2}$-separated,
the $L^{2}$-balls 
\[
B\left(B_{i},\frac{\delta\mathrm{Vol}(X)^{1/2}}{2}\right),\qquad1\leq i\leq N,
\]
 are pairwise disjoint. Moreover, since $\left\Vert B_{i}\right\Vert _{2}\leq\mathrm{Vol}\left(X\right)^{1/2}$,
each of these balls is contained in 
\[
B_{\mathscr{P}(T)}\left(0,\mathrm{Vol}(X)^{1/2}+\frac{\delta\mathrm{Vol}(X)^{1/2}}{2}\right).
\]
The real dimension of $\mathscr{P}\left(T\right)$ is $2\mathrm{dim}_{\mathbb{C}}\mathscr{P}\left(T\right)\leq2m\left(T\right)$.
Comparing the volumes of the disjoint balls with that of the containing
ball gives 
\[
N\leq\left(\frac{\mathrm{Vol}\left(X\right)^{1/2}+\delta\mathrm{Vol}\left(X\right)^{1/2}/2}{\delta\mathrm{Vol}\left(X\right)^{1/2}/2}\right)^{2\mathrm{dim}_{\mathbb{C}}\mathscr{P}\left(T\right)}=\left(1+\frac{2}{\delta}\right)^{2\mathrm{dim}_{\mathbb{C}}\mathscr{P}\left(T\right)}\leq\left(1+\frac{2}{\delta}\right)^{2m\left(T\right)}.
\]
This proves the lemma.
\end{proof}
Since every $B_{i}\in\mathscr{A}(T)$, we may choose a mean-zero observable
$A_{i}$ satisfying $\left\Vert A_{i}\right\Vert _{\infty}\leq1$
and $P_{T}A_{i}=B_{i}$. Thus, for every mean-zero $A$ with $\left\Vert A\right\Vert _{\infty}\leq1$,
there exists $A_{i}$ such that 
\[
\left\Vert P_{T}\left(A-A_{i}\right)\right\Vert _{2}\leq\delta\mathrm{Vol}\left(X\right)^{1/2}.
\]
Applying the estimate in Lemma \ref{lem:lowcontribution inhighT and continuity},
we have 
\begin{equation}
\left|\mathcal{V}_{\Gamma,k}\left(U,A\right)^{1/2}-\mathcal{V}_{\Gamma,k}\left(U,A_{i}\right)^{1/2}\right|\lesssim_{k,r_{0}}\left(\delta^{2}+\eta_{k}\left(T\right)+\mathrm{GeomErr}_{X}\left(R\right)\right)^{1/2}.\label{eq:net bound}
\end{equation}

\subsection{Uniform concentration over the finite net}

Let 
\[
\mathcal{N}_{n}\left(T_{n},\delta_{n}\right)=\left\{ P_{n,T_{n}}A_{n,1},\ldots,P_{n,T_{n}}A_{n,N_{n}}\right\} 
\]
be the finite net constructed in the previous subsection, where 
\[
\left\Vert A_{n,i}\right\Vert _{\infty}\leq1
\]
for every $1\leq i\leq N_{n}$. Then by Lemma \ref{lem:net},
\begin{equation}
N_{n}\leq\left(1+\frac{2}{\delta_{n}}\right)^{2m_{n}(T_{n})},\label{eq:net in n}
\end{equation}
where $m_{n}\left(T_{n}\right)$ denotes the number, counted with
multiplicity, of non-zero eigenvalues of the Laplacian satisfying
$\lambda\leq1/4+T_{n}^{2}$ on $\Gamma_{n}\backslash\mathbb{D}$.
Recall that $m_{A_{n,i}}$ is the median of the function 
\[
U\mapsto\mathcal{V}_{n,k}\left(U,A_{n,i}\right)
\]
on $\mathrm{SU}\left(d_{n,k}\right)$. Since $\mathcal{V}_{n,k}\left(U,A_{n,i}\right)\geq0$,
we have 
\[
m_{A_{n,i}}\leq2\mathbb{E}_{U}\left[\mathcal{V}_{n,k}\left(U,A_{n,i}\right)\right].
\]
By Corollary \ref{cor:EVO(1/d)}, there exists a constant $C_{1}=C_{1}\left(k,r_{0}\right)>0,$
independent of $n$ and $i$, such that 
\begin{equation}
m_{A_{n,i}}\leq C_{1}\left(\frac{1}{d_{n,k}}+\rho_{n}^{2}\right).\label{eq:uniform upper bound for median}
\end{equation}
By (\ref{eq:Gromov-Milman}) and (\ref{eq:Contain1}), we have for
each $i$
\[
\mu_{n,k}\left\{ U\in\mathrm{SU}\left(d_{n,k}\right):\mathcal{V}_{n,k}\left(U,A_{n,i}\right)\leq m_{A_{n,i}}+\frac{4}{\sqrt{d_{n,k}}}\varepsilon_{n}\right\} \geq1-\sqrt{2}\exp\left(-\frac{d_{n,k}\varepsilon_{n}^{2}}{4}\right).
\]
The uniform upper bound for the medians in (\ref{eq:uniform upper bound for median})
implies that if for each $i$, $U$ satisfies 
\[
\mathcal{V}_{n,k}\left(U,A_{n,i}\right)\leq m_{A_{n,i}}+\frac{4}{\sqrt{d_{n,k}}}\varepsilon_{n},
\]
 then for all $i$
\begin{equation}
\mathcal{V}_{n,k}\left(U,A_{n,i}\right)\leq C_{1}\left(\frac{1}{d_{n,k}}+\rho_{n}^{2}\right)+\frac{4\varepsilon_{n}}{\sqrt{d_{n,k}}}.\label{eq:upperbound net Omega}
\end{equation}
Consequently,
\[
\mu_{n,k}\left\{ U\in\mathrm{SU}\left(d_{n,k}\right):\mathcal{V}_{n,k}\left(U,A_{n,i}\right)\leq C_{1}\left(\frac{1}{d_{n,k}}+\rho_{n}^{2}\right)+\frac{4\varepsilon_{n}}{\sqrt{d_{n,k}}}\right\} \geq1-\sqrt{2}\exp\left(-\frac{d_{n,k}\varepsilon_{n}^{2}}{4}\right).
\]
We now apply this estimate simultaneously to all observables $A_{n,i}$.
Define 
\[
\Omega_{n}=\bigcap_{i=1}^{N_{n}}\left\{ U\in\mathrm{SU}\left(d_{n,k}\right):\mathcal{V}_{n,k}\left(U,A_{n,i}\right)\leq C_{1}\left(\frac{1}{d_{n,k}}+\rho_{n}^{2}\right)+\frac{4\varepsilon_{n}}{\sqrt{d_{n,k}}}\right\} .
\]
By the union bound,
\begin{align*}
1-\mu_{n,k}\left(\Omega_{n}\right) & \leq\sum_{i=1}^{N_{n}}\mu_{n,k}\left\{ U\in\mathrm{SU}\left(d_{n,k}\right):\mathcal{V}_{n,k}\left(U,A_{n,i}\right)>C_{1}\left(\frac{1}{d_{n,k}}+\rho_{n}^{2}\right)+\frac{4\varepsilon_{n}}{\sqrt{d_{n,k}}}\right\} \\
 & \leq\sqrt{2}N_{n}\exp\left(-\frac{d_{n,k}\varepsilon_{n}^{2}}{4}\right).
\end{align*}
Using the estimate for the cardinality of the net, we obtain 
\begin{equation}
\mu_{n,k}\left(\Omega_{n}\right)\geq1-\sqrt{2}\left(1+\frac{2}{\delta_{n}}\right)^{2m_{n}(T_{n})}\exp\left(-\frac{d_{n,k}\varepsilon_{n}^{2}}{4}\right).\label{eq:simultaneous bound for i}
\end{equation}

\subsection{Extension to all observables from the finite net}

Let $U\in\Omega_{n}$ , and let $A_{n}$ be any mean-zero observable
satisfying $\left\Vert A_{n}\right\Vert _{\infty}\leq1.$ Then Lemma
\ref{lem:net} implies that there exists $A_{n,i}$ with $\left\Vert A_{n,i}\right\Vert _{\infty}\leq1$
such that 
\[
\left\Vert P_{n,T_{n}}\left(A_{n}-A_{n,i}\right)\right\Vert _{2}\leq\delta_{n}\mathrm{Vol}\left(X_{n}\right)^{1/2}.
\]
By (\ref{eq:net bound}), there exists $C_{2}=C_{2}(k,r_{0})$ such
that 
\[
\left|\mathcal{V}_{n,k}\left(U,A_{n}\right)^{1/2}-\mathcal{V}_{n,k}\left(U,A_{n,i}\right)^{1/2}\right|\leq\left[C_{2}\left(\delta_{n}^{2}+\eta_{k}\left(T_{n}\right)+\mathrm{GeomErr}_{X_{n}}\left(R\right)\right)\right]^{1/2}.
\]
Using the upper bound for $\mathcal{V}_{n,k}\left(U,A_{n,i}\right)$
in (\ref{eq:upperbound net Omega}), for $U\in\Omega_{n}$, 
\[
\mathcal{V}_{n,k}\left(U,A_{n}\right)^{1/2}\leq\left[C_{1}\left(\frac{1}{d_{n,k}}+\rho_{n}^{2}\right)+\frac{4\varepsilon_{n}}{\sqrt{d_{n,k}}}\right]^{1/2}+\left[C_{2}\left(\delta_{n}^{2}+\eta_{k}\left(T_{n}\right)+\mathrm{GeomErr}_{X_{n}}\left(R\right)\right)\right]^{1/2}.
\]
Then, by squaring both sides, we have 
\begin{align}
\mathcal{V}_{n,k}\left(U,A_{n}\right) & \leq\left(\left[C_{1}\left(\frac{1}{d_{n,k}}+\rho_{n}^{2}\right)+\frac{4\varepsilon_{n}}{\sqrt{d_{n,k}}}\right]^{1/2}+\left[C_{2}\left(\delta_{n}^{2}+\eta_{k}\left(T_{n}\right)+\mathrm{GeomErr}_{X_{n}}\left(R\right)\right)\right]^{1/2}\right)^{2}\label{eq:final bound}\\
 & \leq2C_{1}\left(\frac{1}{d_{n,k}}+\rho_{n}^{2}\right)+\frac{8\varepsilon_{n}}{\sqrt{d_{n,k}}}+2C_{2}\left(\delta_{n}^{2}+\eta_{k}\left(T_{n}\right)+\mathrm{GeomErr}_{X_{n}}\left(R\right)\right).\nonumber 
\end{align}
Set 
\[
\tau_{n}=\inf_{R>0}\mathrm{GeomErr}_{X_{n}}\left(R\right).
\]
By Benjamini-Schramm convergence, $\tau_{n}\to0$. This proves Proposition
\ref{prop:main thm}.

\section*{Statement on the Use of Artificial Intelligence}

During the preparation of this work, the author used ChatGPT (GPT-5.6
Sol) to improve language readability and exposition of the manuscript
and for mathematical discussion concerning the calculation of 
\[
\int_{0}^{\infty}\int_{0}^{\pi}\left(\cosh\frac{r}{2}\right)^{-2k}\left(\cosh r-\sinh r\cos2\theta\right)^{-\frac{1+s}{2}}\sinh rd\theta dr
\]
 and Proposition \ref{prop:general meanQV}.

\end{document}